%
%



\def\LaTeX{%
  \let\Begin\begin
  \let\End\end
  \def\Bcenter{\Begin{center}}
  \def\Ecenter{\End{center}}
  \let\Label\label
  \let\salta\relax
  \let\finqui\relax
  \let\futuro\relax}

\def\UK{\def\our{our}\let\sz s}
\def\USA{\def\our{or}\let\sz z}



\LaTeX

\USA


\salta

\documentclass[twoside,12pt]{article}
\setlength{\textheight}{24cm}
\setlength{\textwidth}{16cm}
\setlength{\oddsidemargin}{2mm}
\setlength{\evensidemargin}{2mm}
\setlength{\topmargin}{-15mm}
\parskip2mm


%
%
\usepackage{cite}

\usepackage{color}
\usepackage{amsmath}
\usepackage{amsthm}
\usepackage{amssymb}

\usepackage{amsfonts}
\usepackage{mathrsfs}

\usepackage{hyperref}
\usepackage[mathcal]{euscript}

\usepackage[ulem=normalem,draft]{changes}




%
\newtheorem{theorem}{Theorem}[section]

\newtheorem{corollary}[theorem]{Corollary}

\newtheorem{lemma}[theorem]{Lemma}

\finqui

\def\Beq{\Begin{equation}}
\def\Eeq{\End{equation}}
\def\Bsist{\Begin{eqnarray}}
\def\Esist{\End{eqnarray}}

\def\Bthm{\Begin{theorem}}
\def\Ethm{\End{theorem}}

\def\Brem{\Begin{remark}\rm}
\def\Erem{\End{remark}}

\def\Bdim{\Begin{proof}}
\def\Edim{\End{proof}}
\let\non\nonumber




\def\step #1 \par{\medskip\noindent{\bf #1.}\quad}


\def\lhs{left-hand side}
\def\rhs{right-hand side}



\def\multibold #1{\def\arg{#1}%
  \ifx\arg\pto \let\next\relax
  \else
  \def\next{\expandafter
    \def\csname #1#1#1\endcsname{{\bf #1}}%
    \multibold}%
  \fi \next}

\def\pto{.}

\def\multical #1{\def\arg{#1}%
  \ifx\arg\pto \let\next\relax
  \else
  \def\next{\expandafter
    \def\csname cal#1\endcsname{{\cal #1}}%
    \multical}%
  \fi \next}


\def\multimathop #1 {\def\arg{#1}%
  \ifx\arg\pto \let\next\relax
  \else
  \def\next{\expandafter
    \def\csname #1\endcsname{\mathop{\rm #1}\nolimits}%
    \multimathop}%
  \fi \next}

\multibold
qwertyuiopasdfghjklzxcvbnmQWERTYUIOPASDFGHJKLZXCVBNM.

\multical
QWERTYUIOPASDFGHJKLZXCVBNM.

\multimathop
dist div dom meas sign supp .


\def\Accorpa #1#2 #3 {\gdef #1{\eqref{#2}--\eqref{#3}}%
  \wlog{}\wlog{\string #1 -> #2 - #3}\wlog{}}


\def\graffe #1{\mathopen\{#1\mathclose\}}

\def\<#1>{\mathopen\langle #1\mathclose\rangle}
\def\norma #1{\mathopen \| #1\mathclose \|}


\def\separa{\noalign{\allowbreak}}

\def\iot {\int_0^t}
\def\ioT {\int_0^T}
\def\iO{\int_\Omega}

\def\iQ{\iint_Q}

\def\dt{\partial_t}
\def\dn{\partial_{\bf n}}

\def\checkmmode #1{\relax\ifmmode\hbox{#1}\else{#1}\fi}

\def\aat{\checkmmode{for a.a.~$t\in(0,T)$}}


\def\erre{{\mathbb{R}}}
\def\enne{{\mathbb{N}}}




\def\genspazio #1#2#3#4#5{#1^{#2}(#5,#4;#3)}
\def\spazio #1#2#3{\genspazio {#1}{#2}{#3}T0}

\def\L {\spazio L}
\def\H {\spazio H}

\def\C #1#2{C^{#1}([0,T];#2)}


\def\Lx #1{L^{#1}(\Omega)}
\def\Hx #1{H^{#1}(\Omega)}

\def\Luno{\Lx 1}
\def\Ldue{\Lx 2}
\def\Linfty{\Lx\infty}

\def\Huno{\Hx 1}
\def\Hdue{\Hx 2}

\def\LiQ{L^\infty(Q)}



\let\theta\vartheta

\let\phi\varphi

\let\TeXchi\chi                         
\newbox\chibox
\setbox0 \hbox{\mathsurround0pt $\TeXchi$}
\setbox\chibox \hbox{\raise\dp0 \box 0 }
\def\chi{\copy\chibox}



\def\CX{{\cal X}}

\def\CS{{\cal S}}

\def\CY{{\cal Y}}
\def\CU{{\cal U}}

\def\Uad{{\cal U}_{\rm ad}}
\def\CZ{{\cal Z}}
\def\CN{{\cal N}}
\def\CL{{\cal L}}

\def\us{u^*}
\def\phis{\phi^*}
\def\mus{\mu^*}
\def\ws{w^*}
\def\qs{q^*}
\def\ps{p^*}
\def\rs{r^*}

\def\yh{y^h}
\def\zh{z^h}
\def\zeh{\zeta^h}

\def\xih{\xi^h}
\def\xik{\xi^k}
\def\vh{v^h}
\def\oxih{\overline\xi^h}
\def\oetah{\overline\eta^h}
\def\ovh{\overline v^h}

\def\VD{V^*}
\def\phin{\phi_n}
\def\mun{\mu_n}
\def\wn{w_n}
\def\phiz{\phi_0}


\normalfont

\Begin{document}


\title{{\bf Second-order optimality conditions for the \\sparse optimal control of nonviscous\\ Cahn--Hilliard systems }}
\author{}
\date{}
\maketitle

\Bcenter
\vskip-1.3cm
{\large\bf Pierluigi Colli$^{(1)}$}\\
{\normalsize e-mail: {\tt pierluigi.colli@unipv.it}}\\[.4cm]
{\large\bf J\"urgen Sprekels$^{(2)}$}\\
{\normalsize e-mail: {\tt juergen.sprekels@wias-berlin.de}}\\[.6cm]

$^{(1)}$
{\small Dipartimento di Matematica ``F. Casorati'', Universit\`a di Pavia}\\
{\small and Research Associate at the IMATI -- C.N.R. Pavia}\\ 
{\small via Ferrata 5, 27100 Pavia, Italy}\\[.2cm]
$^{(2)}$
{\small Weierstrass Institute for Applied Analysis and Stochastics}\\
{\small Mohrenstra\ss e 39, 10117 Berlin, Germany}\\[.8cm]

\Ecenter

{
\Begin{abstract}\noindent
In this paper we study the optimal control of an initial-boundary value problem
for the classical nonviscous Cahn--Hilliard system with zero Neumann boundary
conditions. Phase field systems of this type govern the evolution of
diffusive phase transition processes with conserved order parameter. For such systems, optimal control problems have been studied in the past.  
We focus here on the situation when the cost functional of the optimal control problem contains a sparsity-enhancing nondifferentiable term 
like the $L^1$-norm. 
For such cases, we establish first-order necessary and second-order sufficient optimality conditions for locally optimal controls, where in 
the approach to second-order sufficient conditions we employ a technique introduced 
by E.~Casas, C.~Ryll and F.~Tr\"{o}ltzsch  in the paper [{\em SIAM J. Control Optim.} {\bf  53} (2015), 2168--2202]. The main novelty of this
paper is that this method, which has recently been successfully applied to systems of viscous Cahn--Hilliard type,
can be adapted also to the classical nonviscous case. Since in the case without viscosity the solutions to 
the state and adjoint systems turn out to be considerably less regular than in the viscous case, numerous additional technical difficulties 
have to be overcome, and additional conditions have to be imposed. In particular, we have to restrict ourselves to the case when the 
nonlinearity driving the phase 
separation is regular, while in the presence of a viscosity term also nonlinearities of logarithmic type turn could be admitted. 
In addition, the implicit function theorem, which was employed to establish the needed differentiability properties of the 
control-to-state operator in the viscous case, does not apply in our situation and has to be substituted by other arguments.       
\\[2mm]
\normalfont
\\[2mm]

\End{abstract}
}
\salta

\pagestyle{myheadings}
\newcommand\testopari{\sc Colli \ --- \ Sprekels }
\newcommand\testodispari{\sc Optimality conditions with sparsity for the Cahn--Hilliard system}
\markboth{\testodispari}{\testopari} 


\finqui


\section{Introduction}
\label{Intro}
\setcounter{equation}{0}

Let $\Omega\subset \erre^3$ denote some bounded and connected open set with smooth boundary $ \Gamma=\partial\Omega$ (a 
compact hypersurface of class $C^2$), unit outward 
normal ${\bf n}$, and associated outward normal derivative $\,\dn$.
Moreover, let $T>0$ denote some final time, and set
\begin{align*}
&Q_t:=\Omega \times (0,t),\quad \Sigma_t:=\Gamma\times (0,t), \quad \mbox{for }\,t\in(0,T],\quad\mbox{and}\quad
Q:=Q_T,\quad \Sigma:=\Sigma_T.
\end{align*}
We then study the following optimal control problem:

\vspace*{2mm}\noindent
{\bf (CP)} \,\,Minimize the cost functional
\begin{align}
\label{cost}
{\cal J}(\phi, u)\,&:=\,\frac{b_1}2\iQ|\phi -\phi_Q|^2\,+\,\frac{b_2}2 \iO|\phi(T)-\phi_{\Omega}|^2\,+\,\frac{b_3}2\iQ|u|^2
\,+\,\kappa\,G(u)\,,\nonumber\\
&=: J(\phi,u)\,+\,\kappa\,G(u)
\end{align}
subject to the initial-boundary value system 
\begin{align}
\label{ss1}
&\langle \dt\phi,v \rangle\,+\iO\nabla\mu\cdot\nabla v=0 &&\mbox{for all $v\in \Huno$ and a.e. in $(0,T)$}, \\
\label{ss2}
&-\Delta \phi + f'(\phi) = \mu + w  &&\mbox{a.e. in }\,Q,\\
\label{ss3}
&\gamma \dt w + w = u   &&\mbox{a.e. in }\,Q,\\
\label{ss4}
&\dn\phi = 0 &&\mbox{a.e. on }\,\Sigma,\\
\label{ss5}
&\phi(0)=\phi_0,\,\quad w(0)= w_0 &&\mbox{a.e. in }\,\Omega,
\end{align}
\Accorpa\State ss1 ss5
and to the control constraint
\begin{equation}
\label{defUad}
\Uad=\{u\in\CU: \ \underline u(x,t)\le u(x,t)\le \overline u(x,t) \,\mbox{ for a.a. $(x,t)$\, in }\,Q\},
\end{equation}
where the control space is specified by
\begin{equation}
\label{defU}
\CU=L^2(0,T;\Ldue).
\end{equation}
The given bounds $\,\underline u,\overline u\in {\LiQ}\,$ satisfy $\,\underline u\le\overline u$ 
almost everywhere in $Q$. Moreover, the targets $\phi_Q, \, \phi_{\Omega}$ are given functions, $b_1\ge 0$, $b_2\ge 0$, $b_3>0$ are constants,
and
$\kappa>0$ is a constant representing the sparsity parameter. The sparsity-enhancing functional $\,G:L^2(Q)\to\erre\,$ is
nonnegative, continuous and convex. Typically, $G$ has a nondifferentiable form like, e.g.,
\begin{equation}
\label{defg}
G(u)=\|u\|_{L^1(Q)}=\iint_Q|u|\,. 
\end{equation}

The state equations \eqref{ss1}--\eqref{ss2} constitute a classical nonviscous Cahn--Hilliard system in which a number of 
physical constants have been normalized to unity. Notice that \eqref{ss1} is just the weak form of the partial differential equation
$\,\dt\phi-\Delta\mu=0$, where, throughout this paper, $\,\langle\,\cdot\,,\,\cdot\,\rangle\,$ denotes the dual 
pairing between $\Huno$ and its dual $\Huno^*$: actually, in~\eqref{ss1} also the boundary condition $\dn\mu = 0$ is taken into account.
The state variables 
$\,\phi\,$ and $\,\mu\,$  are monitored through 
the input variable $\,w$, which is in turn determined by the action of the control $\,u\,$ via the linear control equation~\eqref{ss3}. 
Equation~\eqref{ss3} models how the ``forcing'' $w$ is generated by the external control $u$. We remark that \eqref{ss4}
could be replaced by much more general partial differential equations modeling the relation between an $L^2$-control $u$ and a forcing $w$.
In the system \eqref{ss1}--\eqref{ss5}, $\,\varphi\,$ plays the role of an \emph{order parameter}, while $\,\mu\,$ is the 
associated \emph{chemical potential}.  
Moreover,  $\gamma $  is a given (uniformly) positive function defined on $\Omega$, and
$ \phi_0$ and $ w_0$ are given initial data. The nonlinearity  $\,f\,$ represents a smooth double-well potential whose
derivative defines the local part of the thermodynamic force driving the evolution of the system. 
A typical form of $\,f\,$ is 
\begin{equation}
\label{deff}
f(\phi)=\frac 14(\phi^2-1)^2.
\end{equation}  

Starting with the seminal paper \cite{EZ}, there exists an abundant literature on the well-posedness and asymptotic behavior of the 
nonviscous Cahn--Hilliard system with zero Neumann and with dynamic boundary conditions that cannot be cited here. 
In spite of this large amount of related literature, we have chosen
to provide a detailed well-posedness analysis of the state system \State, both for the readers' convenience and the fact that 
the system \State\ was apparently not studied before in this particular form in which the control contributes to the chemical potential
through the quantity $w$. 
   
There are also contributions devoted to the optimal control 
of Cahn--Hilliard type systems in various contexts. 
Without claiming to be exhaustive and complete, we mention now some related papers. First, let us
refer to \cite{CGSCC,Duan,HintWeg,Z,ZW} and, in the framework of diffusive models of tumor growth, 
to \cite{CGRS,CSig,CSS1,CSS2,EK1,EK2,GLR,You}. Problems with dynamical boundary conditions were studied in  
\cite{CFGS1,CFGS2,CGSANA,CGSAMO,CGSAnnali,CGSSIAM,CGSconv,CSig,GS}, and convective Cahn--Hilliard systems have been the subject of
\cite{CGSAnnali,CGSSIAM,GS,RoSp,ZL1,ZL2}. In addition, quite a number of works have been dedicated to the study of cases in which 
the Cahn--Hilliard system is coupled to other systems; in this connection, we quote Cahn--Hilliard--Navier--Stokes models
(see \cite{FGS,HK,HKW,HW2,Medjo,Zhao}) and the Cahn--Hilliard--Oono (see \cite{CGRS2,GiRoSi}),  Cahn--Hilliard--Darcy 
(see\cite{ACGW,SpWu}), Cahn--Hilliard--Brinkman (see~\cite{EK1}) and Cahn--Hilliard with curvature effects (see~\cite{CGSS6}) systems.

None of the papers cited above is concerned with the aspect of {\em sparsity}, i.e., the 
possibility that any locally optimal control may vanish in subregions of positive measure of the space-time cylinder $Q$
 that are controlled by the
sparsity parameter $\kappa$. 
In this paper, we focus on sparsity, where we restrict ourselves to the case of {\em full sparsity}
which is connected to the $L^1(Q)$-norm functional $G$ introduced in~\eqref{defg}. Other types of sparsity such as 
{\em directional sparsity with respect to time} 
and {\em directional sparsity with respect to space} (see, e.g.,~\cite{SpTr1}) are not treated in this paper.

Sparsity in the optimal control theory for partial differential equations has become an active field of research. Beginning with 
\cite{stadler2009}, many results on sparse optimal controls for PDEs 
were published. We mention only very few of them with closer relation to our paper, in particular 
\cite{casas_herzog_wachsmuth2017,herzog_obermeier_wachsmuth2015,herzog_stadler_wachsmuth2012} on directional sparsity
 and \cite{casas_troeltzsch2012} on a general theorem for second-order conditions. 
Moreover, we refer to some new trends in the investigation of sparsity, namely, infinite horizon sparse 
optimal control (see, e.g., \cite{Kalise_Kunisch_Rao2017,Kalise_Kunisch_Rao2020}) and fractional order optimal 
control (cf.~\cite{Otarola2020,Otarola_Salgado2018}).

The abovementioned papers concentrated on the first-order optimality conditions for sparse optimal controls of single elliptic and parabolic equations. 
In  \cite{casas_ryll_troeltzsch2013,casas_ryll_troeltzsch2015}, first- and second-order optimality conditions 
have been discussed in the context of sparsity for the (semilinear) system of  FitzHugh--Nagumo equations. 
More recently, sparsity of optimal controls for reaction-diffusion systems of Cahn--Hilliard type have been addressed in 
\cite{CSS4,Garcke_etal2021,SpTr1}.  Moreover, we refer to the measure control of the Navier--Stokes system 
studied in \cite{Casas_Kunisch2021}.
Second-order sufficient optimality conditions for sparse controls 
for the viscous Cahn--Hilliard system were recently addressed in 
\cite{CoSpTr}. 

Second-order sufficient optimality 
conditions are based on a condition of coercivity that is required to hold for the smooth part  $\,J\,$ 
of $\,{\cal J}\,$ in a certain {\em critical cone}. The nonsmooth part $\,G\,$ contributes to sufficiency by its convexity. 
For the strength of sufficient conditions it is crucial that the critical cone be as small as possible. In their paper
\cite{casas_ryll_troeltzsch2015}, Casas--Ryll--Tr\"oltzsch devised a technique by means of which a very advantageous (i.e., small)
critical cone can be chosen. This method was originally introduced  for a class of semilinear second-order parabolic problems
with smooth nonlinearities. In the recent papers \cite{SpTr2,SpTr3,CoSpTr} it has been demonstrated that 
it can be adapted correspondingly to the sparse optimal control of Allen--Cahn systems
with dynamic boundary conditions \cite{SpTr2}, to a large class of systems modeling tumor growth \cite{SpTr3}, and to the viscous Cahn--Hilliard system
\cite{CoSpTr},
where in all of these papers logarithmic nonlinearities could be admitted.

It is the main aim and novelty of this work to show that also the classical nonviscous Cahn--Hilliard system can be treated
accordingly. This is by no means obvious, since in the nonviscous case the solutions to the state and the adjoint state systems 
enjoy less regularity than in the viscous one. And indeed, it turns out that the needed analytic effort is considerably
larger than in the viscous case, while the optimization part of the argument changes only little.  

In particular, one of the key elements of the technique is to show that the reduced cost functional associated with 
the smooth part $\,J\,$ of $\,{\cal J}\,$ is twice
continuously differentiable, which in turn requires that the  
control-to-state operator is twice continuously Fr\'echet differentiable between $\,\CU=\L2{L^2(\Omega)}\,$ and a suitable Banach space.
For this to be the case, it seems mandatory that the phase variable $\,\phi\,$ satisfies a  
so-called {\em uniform separation property}. While such a condition can in the viscous case be shown also for logarithmic potentials, 
a corresponding result seems in the nonviscous case to be available only for regular potentials. We therefore have to restrict
our analysis to such nonlinearities in this paper, thereby excluding logarithmic potentials.
 
Another difficulty arises from the fact that the first component of the solution triple to the linearized
state system (see \eqref{ls1}--\eqref{ls5} below) is not known to be bounded. This entails that the technique employed 
in \cite{CoSpTr}, which is based on the good differentiability 
properties of Nemytskii operators on $L^\infty$ spaces and makes use of the implicit function theorem,
does not carry over to the nonviscous case. All this has the consequence that the proof of the needed twice continuous Fr\'echet
differentiability of the control-to-state operator is technically quite challenging and sometimes a bit tedious. 
Nevertheless, it turns out
that all difficulties can be overcome.  

The paper is organized as follows. In the following section, we formulate the general assumptions and study the state
system, proving the existence of a unique solution. We also show the local Lipschitz continuity of the control-to-state operator. 
In Section 3, we then prove that  the control-to-state operator is twice continuously Fr\'echet differentiable between appropriate 
Banach spaces. The proofs in this section require the main analytical effort of this paper. In Section 4, 
we investigate the control problem {\bf (CP)} with sparsity. Besides analyzing the associated adjoint problem, we
derive the first-order necessary optimality conditions. The final section then brings the derivation of the announced second-order 
sufficient optimality conditions for controls that are locally optimal in the sense of $L^2(Q)$.

Prior to this, let us fix some notation.
For any Banach space $X$, we denote by \,$\|\,\cdot\,\|_X$ and $X^*$   
the corresponding norm and its dual space, respectively. 
For two Banach spaces $X$ and $Y$ that are both continuously embedded in some topological vector space~$Z$, we consider the linear space
$X\cap Y$ that becomes a Banach space if equipped with its natural norm $\norma v_{X\cap Y}:=\norma v_X+\norma v_Y\,$  for $v\in X\cap Y$.
The standard Lebesgue and Sobolev spaces defined on $\Omega$ are, for
$1\le p\le\infty$ and $\,m \in \enne \cup \{0\}$, denoted by $L^p(\Omega)$ and $W^{m,p}(\Omega)$, respectively. 
If $p=2$, we also use the usual notation $H^m(\Omega):= W^{m,2}(\Omega)$. Moreover, for convenience, we denote the norm of $\,L^p(\Omega)\,$ by $\,\|\,\cdot\,\|_p\,$ for $1\le p\le\infty$, and we set
\begin{align*}
  & H := \Ldue , \quad V := \Huno, \quad W:=\bigl\{v\in \Hdue : \ \dn v =0 \, \hbox{ on }\, \Gamma \bigr\},
  \end{align*}
where we denote by $(\,\cdot\,,\,\cdot\,)$ the natural inner product in $\,H$.
As usual, $H$ is identified with a subspace of the dual space $\VD$ according to the identity
\begin{align*}
	\langle u,v\rangle =(u,v)
	\quad\mbox{for every $u\in H$ and $v\in V$}.
	\end{align*}
We then have the Hilbert triple $V \subset H \subset \VD$ with dense and compact embeddings.

We close this section by introducing a convention concerning the constants used in estimates within this paper: we denote by $\,C\,$ any 
positive constant that depends only on the given data occurring in the state system and in the cost functional, as well as 
on a constant that bounds the $L^2(Q)$--norms of the elements of $\Uad$. The actual value of 
such generic constants $\,C\,$ 	may 
change from formula to formula or even within formulas. Finally, the notation $C_\delta$ indicates a positive constant that
additionally depends on the quantity $\delta$.   


\section{Properties of the state system}
\setcounter{equation}{0}

\subsection{Notation and assumptions}
In this section, we formulate the general assumptions for the data of the state system \State, and we introduce some known tools
for later use. Throughout this paper, we generally assume:
\begin{description}
\item[(A1)] \,\,$f=f_1+f_2$, where $f_1 \in C^5(\erre)$ is a convex and nonnegative function with $f_1(0)=0$ and $f_2\in C^5(\erre)$ has 
a Lipschitz continuous first derivative $f_2'$ on $\erre$.
\item[(A2)] \,\,$\gamma>0\,$ is a constant.
Moreover, $w_0\in V$ and \,$\phi_0 \in H^3(\Omega)\cap W $.
\item[(A3)] \,\,$R>0$  is a fixed constant, and 
$\,\,   \CU_R :=\{u\in L^2(Q): \,\,\|u\|_{L^2(Q)}\,<R\}$.
\end{description} 

\vspace*{1mm} \noindent
From the condition~{\bf (A1)} it follows that $f_1'$ is monotone increasing on $\erre$ and induces a maximal monotone operator 
in $\erre \times \erre$.
Note that {\bf (A1)} is fulfilled if $\,f\,$ is given by the potential \eqref{deff}, and  the 
condition $\phi_0 \in H^3(\Omega)\cap W $ implies that $\phi_0\in C^0(\overline \Omega)$. Moreover, 
the mean value 
\Beq
  m_0 := \frac 1 {|\Omega|}\iO \phiz ,
  \label{meanz}
\Eeq
where $\,|\Omega|\,$ denotes the Lebesgue measure of $\,\Omega$, belongs to a bounded interval in $\erre$. 
In the following, we use the notation $\overline  v$ to denote the mean value of a generic function $v\in\Luno$.
More generally, we set
\Beq
  \overline v : = \frac 1 {|\Omega|} \, \langle v , 1 \rangle
  \quad \hbox{for every $v\in V^*$},
  \label{defmean}
\Eeq
noting that the constant function 1 is an element of $V$. Clearly, $\overline  v$ is the usual mean value of $v$ if $v\in H$.

Next, we recall an important tool which is commonly used when working with problems connected to the Cahn--Hilliard equation.
To this end, consider the weak formulation of the Poisson equation $-\Delta z=\zeta$
with homogeneous Neumann boundary conditions. 
Namely, for a given $\zeta\in\VD$ (which does not necessarily belong to $H$), we consider the problem
of finding
\begin{align}
  z \in V
  \quad \hbox{such that} \quad
  \iO \nabla z \cdot \nabla v
  = \langle \zeta , v \rangle
  \quad \hbox{for every $v\in V$}.
  \label{neumann}
\end{align}
Since $\,\Omega\,$ is connected and smooth, it is well known that the above problem admits solutions $z\in V$ if and only if $\zeta$ has
zero mean value. Hence, we can introduce the following solution operator $\CN$ by setting
\begin{align}
  & \CN: \dom(\CN) := \graffe{\zeta\in\VD:\ \overline \zeta=0} \to \graffe{z\in V:\ \overline z=0},
  \quad 
  \CN: \zeta \mapsto z,
  \label{defN}
\end{align}
where $z$ is the unique solution to \eqref{neumann} satisfying $\overline z=0$.
It turns out that $\CN$ is an isomorphism between the above spaces, and it follows that the formula
\Beq
  \|\zeta\|_*^2 := \iO|\nabla\CN(\zeta-\overline\zeta)|^2 + |\overline \zeta|^2
  \quad \hbox{for every $\zeta\in\VD$}
  \label{normaVp}
\Eeq
defines a Hilbert norm in $\VD$ that is equivalent to the standard dual norm of~$\VD$. In particular, there is a constant $C_\Omega>0$, which depends 
only on $\,\Omega$, such that
\begin{equation}
\label{umba1}
|\langle \zeta,v \rangle|\,\le\,C_\Omega\,\|\zeta\|_*\,\|v\|_V\quad\mbox{for all $\,\zeta\in\VD\,$ and \,$v\in V$}.
\end{equation}
Moreover, from the Young, Poincar\'e and Sobolev inequalities, elliptic estimates, and Ehrling's lemma, we have the estimates
\begin{align}
&ab\,\le\,\frac {\delta} 2|a|^2\,+\,\frac 1{2\delta}|b|^2 \quad\mbox{for all \,$a,b\in\erre$\,and \,$\delta>0$},
\label{Young}
\\[2mm]
& \|v\|_V
  \leq C_\Omega \, \bigl( \|\nabla v\|_{H\times H\times H} + |\overline v| \bigr)
  \quad \hbox{for every $v\in V$,}
  \label{poincare}
  \\[2mm]
  \separa
  & \|v\|_W \leq C_\Omega \bigl( \|\Delta v\|_H + \|v\|_* \bigr) \quad \hbox{for every $v\in W$,}
  \label{elliptic}
    \\[2mm]
  & \|v\|_p
  \leq \delta \, \|\nabla v\|_{H\times H\times H} + C_{\Omega,p,\delta} \, \|v\|_*
  \quad \hbox{for every $v\in V$, $p\in[1,6)$ and $\delta>0$,}
  \label{compactV}
  \\[2mm]
  \separa
& \|v\|_V
  \leq \delta \, \|\Delta v\|_H + C_{\Omega,\delta} \, \|v\|_*
  \quad \hbox{for every $v\in W$ and $\delta>0$}.
  \label{compactW}
\end{align}
In addition, from the above properties there follow the identities
\begin{align}
  & \iO \nabla\CN\zeta \cdot \nabla v
  = \langle \zeta , v \rangle 
  \quad  \hbox{for every $\zeta\in\dom(\CN)$ and $v\in V$,}
  \label{dadefN}
  \\
  & \langle\zeta , \CN\xi \rangle
  = \langle \xi , \CN\zeta \rangle
  \quad \hbox{for every $\zeta,\xi\in\dom(\CN)$,}
  \label{simmN}
  \\
  & \langle \zeta , \CN\zeta \rangle
  = \iO |\nabla\CN\zeta|^2
  = \|\zeta\|_*^2
  \quad \hbox{for every $\zeta\in\dom(\CN)$.}
  \label{danormaVp}
\end{align}
Moreover, we point out the equality
\Beq
  \langle \dt\zeta(t) , \CN\zeta(t) \rangle
  = \frac 12 \, \frac d{dt} \|\zeta(t)\|_*^2
  \quad \aat\,,
  \label{propN} 
\Eeq
which holds true for every $\zeta\in H^1(0,T;\VD)$ satisfying $\overline\zeta=0$ for a.e. $t\in (0,T)$.

\subsection{Existence for the state system}
In this section, we are going to prove an existence result for the state system. Prior to this, we notice
that, thanks to the linear equation~\eqref{ss3} and the second initial condition in \eqref{ss5}, $w\,$ can be explicitly written in terms of $\,u\,$ by means of the variation of constants formula
\Beq
\label{varconst} 
w(x,t)= w_0(x) \exp (-t/{\gamma}) + \int_0^t \exp(-(t-s)/{\gamma} ) u(x,s) ds , \quad \hbox{a.e.}\  (x,t) \in Q.
\Eeq
We have the following result.
\Bthm
\label{Exist}
Suppose that {\bf (A1)}--{\bf (A3)} are fulfilled. Then the state system \State\ has for every $\,u\in\L2H\,$
a unique solution triple $\,(\phi,\mu,w)\,$ satisfying
\begin{align}
\label{regphi}
& \phi \in  W^{1,\infty}(0,T;\VD)\cap H^1(0,T;V) \cap L^\infty(0,T;W)\cap C^0(\overline Q), \\
\label{regmu}
& \mu \in H^1(0,T;\VD)\cap L^\infty(0,T;V) \cap L^2(0,T;W\cap H^3(\Omega)),\\ 
\label{regw}
&w \in \H1 H.
\end{align}
In addition, there is a constant $K_1>0$, which depends only on $\norma{u}_{\L2H}$ and the data of the state system, such that 
\begin{align}
\label{ssbound1}
&\|\phi\|_{ W^{1,\infty}(0,T;\VD) \cap \H1V \cap \L\infty W\cap C^0(\overline Q) } \nonumber\\
&+\,\|\mu\|_{H^1(0,T;\VD)\cap L^\infty(0,T;V)\cap\L 2{W\cap H^3(\Omega)}} \,
+ \, \| w \|_{ \H1 H }       \,\le\,K_1\,  
\end{align} 
whenever $(\phi,\mu,w)$ is the solution to the state system associated with $\,u$. 
\Ethm
\begin{proof}
Although the proof of the above result seems to be pretty standard by now, we carry it out for the reader's 
convenience. We argue by a Faedo--Galerkin approximation. To this end, consider the  
eigenvalues $\,\{\lambda_j\}_{j\in\enne}\,$ of the eigenvalue problem 
$$-\Delta v=\lambda v \quad\mbox{in }\,\Omega, \qquad \dn v=0 \quad\mbox{on }\,\Gamma,$$
and let $\,\{e_j\}_{j\in\enne}\subset W\,$ be the associated eigenfunctions, normalized by $\,\|e_j\|_H=1$. Then 
\begin{align*}
&0=\lambda_1<\lambda_2\le \ldots, \qquad \lim_{j\to\infty}\lambda_j=+\infty,\\
&\iO e_je_k=\iO \nabla e_j\cdot\nabla e_k=0\quad \mbox{for $\,j\not= k$}.
\end{align*}
We then define the $n$-dimensional spaces $\,V_n:={\rm span}\{e_1,\ldots,e_n\}$ for $\,n\in\enne$, where we observe
that $V_1$ is just the space of constant functions on $\,\Omega$. 
It is well known that the union of these spaces is dense in both $\,H\,$ and $\,V$. Notice also that 
\begin{equation}
\label{NvinVn}
\CN v\in V_n \quad\mbox{for every $v\in V_n$ with $\,\overline v=0$}.
\end{equation}  
Indeed, if $v\in V_n$ and $\overline v=0$ then $v=\sum_{j=2}^nc_je_j$ with suitable $c_j\in\erre$, $j=1,\ldots n$,
and $z:=\sum_{j=2}^n \lambda_j^{-1}c_je_j\in V_n$ satisfies $\overline z=0$ and $-\Delta z=v$, that is, $z=\CN v$.

The approximating $n$-dimensional problem then reads as follows:
find functions
\begin{equation}
\label{discrete}
\phi_n(x,t)=\sum_{j=1}^n \phi_{nj}(t)e_j(x),\quad \mu_n(x,t)=\sum_{j=1}^n \mu_{nj}(t)e_j(x),\quad w_n(x,t)=\sum_{j=1}^n w_{nj}(t)e_j(x),
\end{equation}
such that 
\begin{align}
\label{ss1n}
&\langle \dt\phin,v \rangle+\iO \nabla\mun\cdot\nabla v=0 &&\mbox{for all $\,v\in V_n\,$, a.e. 
in $(0,T),$}\\
\label{ss2n}
&\iO \nabla\phin\cdot\nabla v+(f'(\phin),v)-(\wn,v) = (\mun,v) 
&&\mbox{for all $\,v\in V_n\,$, a.e. in $(0,T),$}\\
\label{ss3n}
&\gamma \,(\dt\wn,v) +(\wn,v)=(u,v)&&\mbox{for all $\,v\in V_n\,$, a.e. in $(0,T),$}\\
\label{ss4n}
&\phin(0)=P_n(\phiz), \quad \wn(0)=P_n(w_0),&&\mbox{a.e. in $\,\Omega$,}
\end{align}
where $P_n$ denotes the $H$-orthogonal projection onto $V_n$. Then $P_n(v)=\sum_{j=1}^n (v,e_j)e_j$ for every $v\in H$, and we have  
(see \cite[formula~(3.14)]{CGSS6}), with a constant $C_\Omega>0$ depending only on $\,\Omega$,
\begin{equation}
\label{Fiete}
\|P_n (v)\|_Y\,\le\,C_\Omega \|v\|_Y \quad\mbox{for every $\,v\in Y$, where }\,Y\in\{H,V,W\}.
\end{equation}                                             

Next, we  insert $v=e_k$ in all of the equations \eqref{ss1n}--\eqref{ss4n}, for $k=1,\ldots,n$, obtaining the system
\begin{align}
\label{ODE1}
&\frac{d}{dt} \phi_{nk}+\lambda_k\,\mu_{nk}=0 &&\mbox{a.e. in }\,(0,T),\\
\label{ODE2}
&\mu_{nk}=\lambda_k\phi_{nk}+(f'(\phin),e_k)-w_{nk} &&\mbox{a.e. in}\,(0,T),\\
\label{ODE3}
&\gamma \,\frac{d}{dt} w_{nk}+w_{nk}=(u,e_k) &&\mbox{a.e. in }\,(0,T),\\
\label{ODE4}
&\phi_{nk}(0)=(\phiz,e_k),\quad w_{nk}(0)=(w_0,e_k).&&{}
\end{align}
Now insert for $\mu_{nk}$ in \eqref{ODE1}, using \eqref{ODE2}. We then arrive at an initial value problem for an explicit ODE system in
the $2n$ unknowns $\,(\phi_{n1},...,\phi_{nn}, w_{n1},...,w_{nn})\,$ with locally Lipschitz continuous nonlinearities and coefficient
functions in $L^2(0,T)$. By Carath\'eodory's theorem, this ODE system has a unique maximal solution 
belonging to $H^1(0,T_n;\erre^{2n})$ for some $T_n\in (0,T]$. This solution in turn 
uniquely determines
via \eqref{ODE2} and \eqref{discrete} a triple $(\phin,\mun,\wn)\in (H^1(0,T_n;V_n))^3 $ 
that
solves \eqref{ss1n}--\eqref{ss4n} on $\Omega\times [0,T_n]$, with the regularity of $\mun $ following from \eqref{ODE2} and {\bf (A1)}.
We show that one can take $T_n=T$. We do this by deriving a series of a priori
estimates for the finite-dimensional approximations. In the following, $C>0$ denotes constants that may depend on \,$\norma{u}_{\L2H}$\, 
and the data of the 
state system, but not on $n\in\enne$.   

\step
First estimate

Testing \eqref{ss3n} by $\,\dt w_n$, with the help of \eqref{Young} we immediately get the estimate
\begin{equation}
\label{esti1}
\|w_n\|_{H^1(0,T_n;H)}\le C\,.
\end{equation}
Then, we test \eqref{ss1n} by $\phin\in V_n$ and \eqref{ss2n} by $\,-\Delta\phin\in V_n$, add, and integrate over $[0,t]$ where $t\in (0,T_n]$.
After a cancellation and reorganisation of terms, we obtain that
\begin{align*}
&\frac 12\|\phin(t)\|_H^2+\iint_{Q_t} |\Delta\phin|^2 +\iint_{Q_t} f_1''(\phin)|\nabla\phin|^2\\
& \le\,\frac 12 \|\phin(0)\|_H^2 -\iint_{Q_t} \wn\,\Delta\phin +\iint_{Q_t} f_2'(\phin)\Delta\phin\,.
\end{align*}
By the convexity of $f_1$, the last term on the \lhs\ is nonnegative. Moreover, owing to {\bf (A1)}, we have $\,|f_2'(\phin)|\le C(1+|\phin|)$ a.e.
in $\Omega\times (0,T_n)$. In view of \eqref{esti1} and Young's inequality, the sum of the last two terms on the \rhs\ is therefore bounded by
$$
\frac 12 \iint_{Q_t}|\Delta\phin|^2+C\int_0^t \left(1+\|\phin(s)\|^2_H\right)\,ds\,.
$$  
Consequently, by Gronwall's lemma, and using the estimate \eqref{elliptic}, we infer that 
\begin{equation}
\label{esti2}
\|\phin\|_{L^\infty(0,T_n;H)\cap L^2(0,T_n;W)}\,\le\,C\,.
\end{equation}
 
We  can draw an important consequence from \eqref{esti1} and \eqref{esti2}: indeed, by  a standard argument it follows from these bounds that the local
solution $(\phi_{n1},...,\phi_{nn},w_{n1},...,w_{nn})$ to the ODE system resulting from \eqref{ODE1}--\eqref{ODE4} is uniformly bounded and thus,
by its maximality, global. Therefore, it  must exist on the whole interval $[0,T]$, that is, we have $T_n=T$. We will exploit this fact in the
following estimates.

\step 
Second estimate

Next, we recall that the constant function $v=1$ belongs to $V_1$. Inserting it in \eqref{ss1n} yields that $\,\overline{\dt\phin}=0\,$ a.e. in
$(0,T)$, which entails that $\,\CN(\dt\phin)\,$ is well defined and belongs to $\,V_n$. We now insert $v=\CN(\dt\phin)$ in \eqref{ss1n} and
$v=\dt\phin$ in \eqref{ss2n}, add, and integrate over $[0,t]$ where $t\in(0,T]$. Using \eqref{dadefN} and \eqref{danormaVp}, and noting the
cancellation of two terms, we obtain the identity
\begin{align*}
&\iot\|\dt\phin(s)\|^2_*\,ds+\frac 12\iO|\nabla\phin(t)|^2+\iO f_1(\phin(t))\\
&=\frac 12\|\nabla P_n(\phiz)\|^2+\iO f_1(P_n(\phiz)) +\iint_{Q_t} w_n\dt\phin -\iint_{Q_t} f_2'(\phin)\dt\phin \,.
\end{align*}

By {\bf (A1)}, the third term on the \lhs\ is nonnegative. Moreover, by \eqref{Fiete}, the first summand on the \rhs\ is bounded, and we
have, using the continuity of the embedding $\Hdue\subset C^0(\overline\Omega)$, that \,$\|P_n(\phiz)\|_{C^0(\overline\Omega)}
\le C\|P_n(\phiz)\|_W\le C\|\phiz\|_W$.
This obviously implies that the sequence \,$\{\iO f_1(P_n(\phiz))\}\,$ is bounded. Moreover, we obtain from \eqref{esti1}, \eqref{esti2}, and 
Young's inequality, that
\begin{align*}
\iint_{Q_t} \wn\dt\phin=\iO \wn(t)\phin(t)-\iO P_n(w_0)P_n(\phiz)-\iint_{Q_t}\phin\dt\wn\,\le\,C\,.
\end{align*}

Finally, we infer from \eqref{umba1}, \eqref{esti2}, and Young's inequality, that
\begin{align*}
&-\iint_{Q_t} f_2'(\phin)\dt\phin\,\le\,C\iot \|\dt\phin(s)\|_*\,\|f_2'(\phin(s))\|_V\,ds\\
&\le \,\frac 12\iot\|\dt\phin(s)\|_*^2\,ds\,+\,C\iint_{Q_t}\left(|f_2'(\phin)|^2\,+\,|f_2''(\phin)|^2\,|\nabla\phin|^2\right)\\
&\le \,\frac 12\iot\|\dt\phin(s)\|_*^2\,ds \,+C\,,
\end{align*}
since, owing to {\bf (A1)}, we have $\,|f_2'(\phin)|\le C(1+|\phin|)\,$ and \,$|f_2''(\phin)|\le C\,$ a.e. in $Q$. Combining the 
above estimates, we have therefore shown that
\begin{equation}
\label{esti3}
\|\phin\|_{H^1(0,T;\VD)\cap L^\infty(0,T;V)}\,+\,\|f_1(\phin)\|_{L^\infty(0,T;L^1(\Omega))}\,\le\,C\,.
\end{equation}

\step
Third estimate

At this point, we recall that $\,\overline{\dt\phin}=0$ a.e. in $(0,T)$, which implies that~(cf.~\eqref{meanz})
\begin{align}
\label{MV}
\overline{\phin(t)}=\overline{P_n(\phiz)}=  \frac1{|\Omega|} \int_\Omega (\phiz, e_1) e_1=
\overline{\phiz}\, \norma{e_1}_H^2 = m_0 
\quad\mbox{for all }\,t\in [0,T].
\end{align}
For almost every $t\in (0,T)$, we now insert $v=\CN(\phin(t)-m_0)$ in \eqref{ss1n} and $v=\phin(t)-m_0$ in \eqref{ss2n}, and add the results.
We obtain
\begin{align}
\label{e31}
&\iO |\nabla(\phin(t)-m_0)|^2\,+\iO f_1'(\phin(t))(\phin(t)-m_0)\non\\
&=\,-\langle \dt\phin(t),\CN(\phin(t)-m_0) \rangle\,+\iO (\wn(t)-f_2'(\phin(t))(\phin(t)-m_0)\non\\
&\le\,C\,\|\dt\phin(t)\|_*\,\|\CN(\phin(t)-m_0)\|_V
\non\\
&\ \quad{}+\,(\|\wn(t)\|_H+\|f_2'(\phin(t))\|_H)\,(\|\phin(t)-m_0\|_H)\,.
\end{align}
Owing to \eqref{esti3} and to the bounds for $\CN(\phin(t)-m_0$ implied by \eqref{esti2} and \eqref{MV}, it follows that 
the first summand on the \rhs\ is bounded in $L^2(0,T)$. Moreover, the second summand is already known to be bounded in $L^\infty(0,T)$.  

Now recall that $\,f_1'\,$ is monotone increasing and that \eqref{MV} holds. Then there exist constants $\,\delta_0>0\,$ and 
$\,C_0>0\,$ such that 
\begin{equation}
\label{Zelik}
f_1'(r)(r-m_0)\ge \delta_0|f_1'(r)|-C_0 \quad\mbox{for every $\,r\in\erre$}.
\end{equation}%
For this estimate we refer to~\cite[Appendix, Prop. A.1]{MiZe}
(see also the detailed proof given in~\cite[p.~908]{GiMiSchi}).
Applying \eqref{Zelik}, we therefore can infer from \eqref{e31} that
\begin{equation}
\label{e32}
\|f_1'(\phin)\|_{L^2(0,T;L^1(\Omega))}\,\le\,C\,.
\end{equation}

Next, we insert the constant function $\,1\in V_1\,$ in \eqref{ss2n}.  We obtain, 
for a.e. $t\in(0,T)$,
\begin{align}
\label{Uwe}
\iO f'_1(\phin(t)) + \iO (f'_2(\phin(t)) - \wn(t)) \,=\, |\Omega|\,\overline{\mun(t)}.
\end{align}
By~\eqref{e32} and previous estimates, both summands on the \lhs\ are bounded in $L^2(0,T)$. Then we conclude that
\begin{equation}
\label{e33}
\|\overline{\mun}\|_{L^2(0,T)}\,\le\,C\,.
\end{equation}

At this point, we test \eqref{ss1n} by $v=\mun(t)-\overline{\mun(t)}$, which has zero mean value. It then follows 
from Young's inequality and \eqref{esti3} that
\begin{align*}
&\iint_Q|\nabla\mun|^2\,=\,-\ioT\langle \dt\phin(t),\mun(t)-\overline{\mun(t)}\rangle\,dt\non\\
&\le\,C\ioT \|\dt\phin(t)\|_*\,\|(\mun-\overline{\mun})(t)\|_V\,dt\,\le\,C\,+\,\frac 12\iint_Q|\nabla\mun|^2\,,
\end{align*}  
whence, invoking also \eqref{e33} and the Poincar\'e 
inequality \eqref{poincare}, we arrive at 
\begin{equation}
\label{esti4}
\|\mun\|_{L^2(0,T;V)}\,\le\,C\,.
\end{equation}

\step
Fourth estimate

We now differentiate both the equations \eqref{ss1n} and \eqref{ss2n} with respect to $\,t$, then test the first of the resulting 
equations by $v= \CN(\dt\phin(t))$ and the second by $v=\dt\phin(t)$. Addition and integration over $(0,t)$, where $t\in (0,T]$, and
use of the properties of the operator $\,\CN$, lead after the cancellation of two terms to the identity
\begin{align}
\label{e41}
&\frac 12 \|\dt\phin(t)\|_*^2\,+\iint_{Q_t} |\nabla\dt\phin|^2\,+\iint_{Q_t} f_1''(\phin)|\dt\phin|^2 \non\\
&=\,\frac 12\|\dt\phin(0)\|_*^2 \,-\iint_{Q_t} f_2''(\phin)|\dt\phin|^2\,+\iint_{Q_t} \dt\wn\,\dt\phin\,.
\end{align} 
By the convexity of $\,f_1$, the third term on the \lhs\ is nonnegative. Moreover, the sum of the last two terms on the \rhs, which
we denote by $\,I$, can be estimated as follows:
\begin{align*}
I\,&\le \,C\,+C\iint_{Q_t}|\dt\phin|^2\,\le\,\frac 12 \iint_{Q_t}|\nabla\dt\phin|^2 \,+\,C\iot\|\dt\phin(s)\|_*^2\,ds \\
&\le \,C\,+\,\frac 12\iint_{Q_t} |\nabla\dt\phin|^2\,.
\end{align*}
Here we have used \eqref{esti1}, \eqref{esti3}, {\bf (A1)}, Young's inequality, and the compactness inequality \eqref{compactV} with $p=2$.

For the initial value we have
\begin{align*}
\|\dt\phin(0)\|_*^2\,&=\,-\iO \mun(0)\,\dt\phin(0)\,=\,-\iO (-\Delta(\phin(0))+f'(\phin(0))-\wn(0))\dt\phin(0)\\
&\le\,C\,\|\dt\phin(0)\|_*\,\|-\Delta(\phin(0))+f'(\phin(0))-\wn(0)\|_V\,.
\end{align*}
We claim that the second factor on the \rhs\ in bounded. Indeed, by \eqref{Fiete}, we have \,$\|\wn(0)\|_V\,=\,\|P_nw_0\|_V\,\le\,C\|w_0\|_V
<+\infty$, \,since (see {\bf (A2)}) $\,w_0\in V$. Moreover, we have already shown above in the second estimate that the sequence $\,\{P_n(\phiz)\}\,$
is bounded in $C^0(\overline\Omega)$. Therefore, as $\,\phiz\in V\,$ and again by \eqref{Fiete},
\begin{align*}
\|f'(\phin(0))\|_V^2\,&=\,\|f'(P_n (\phiz))\|_V^2\,=\,\iO(|f'(P_n(\phiz))|^2\,+\,|f''(P_n(\phiz))|^2\,|\nabla P_n(\phiz)|^2)\\
&\le C\,+\,C\|P_n(\phiz)\|_V^2 \,\le\,C\,+\,C\|\phiz\|_V^2\,<+\infty\,.
\end{align*}
Finally, we have
\begin{align*}
&\Delta P_n(\phiz) \,=\,\sum_{j=1}^n (\phiz,e_j)\Delta e_j \,=\,-\sum_{j=1}^n \lambda_j(\phiz,e_j)e_j\,, \quad\,\,
\nabla\Delta P_n(\phiz)\,=\,-\sum_{j=1}^n \lambda_j(\phiz,e_j)\nabla e_j\,,
\end{align*}
whence
$$
\|\Delta P_n(\phiz)\|^2_V \,=\,\sum_{j=1}^n (\lambda_j^2+\lambda_j^3)\,|(\phiz,e_j)|^2 \,\le\,\sum_{j=1}^\infty 
(\lambda_j^2+\lambda_j^3)\,|(\phiz,e_j)|^2\,<\,+\infty\,,
$$  
since, by {\bf (A2)}, $\,\phiz\in H^3(\Omega)\cap W$. In conclusion, we have shown that
\begin{equation}
\label{esti5}
\|\phin\|_{W^{1,\infty}(0,T;V^*)\cap H^1(0,T;V)}\,\le\,C\,.
\end{equation}
In particular, we now see that the \rhs\ of \eqref{e31} is even bounded in $L^\infty(0,T)$, so that
\begin{equation*}
\|f_1'(\phin)\|_{L^\infty(0,T;L^1(\Omega))}\,\le\,C\,,
\end{equation*} 
and it follows from \eqref{Uwe} that
\begin{equation}
\label{Uwedue}
\|\overline{\mun}\|_{L^\infty(0,T)}\,\le\,C\,.
\end{equation}
At this point, we test \eqref{ss2n} by $v=(\mun-\overline\mun)(t)$, without integrating over time. As at the end of the
third estimate, it then follows from \eqref{Uwedue}
and Poincar\'e's inequality \eqref{poincare} that
\begin{equation}
\label{esti6}
\|\mun\|_{L^\infty(0,T;V)}\,\le\,C\,.
\end{equation}

Next, we test \eqref{ss2n} for a.e. $t\in(0,T)$ by $v=-\Delta\phin(t)\in V_n$, without integrating over time. We obtain
\begin{align*}
\iO|\Delta\phin(t)|^2\,+\iO f_1''(\phin(t))|\nabla\phin(t)|^2\,=\,-\iO(\wn(t)+\mun(t)-f_2'(\phin(t))\,\Delta\phin(t)\,,
\end{align*}
whence, using Young's inequality, \eqref{esti1}, \eqref{esti3}, \eqref{esti6}, and elliptic regularity, 
\begin{equation}
\label{esti7}
\|\phin\|_{L^\infty(0,T;W)}\,\le\,C\,.
\end{equation} 

\step
Existence

By virtue of the uniform estimates shown in the previous steps, there exists a triple $(\phi,\mu,w)$ such that (possibly only on 
a suitable subsequence which is again labeled by $n\in\enne$)
\begin{align}
\label{conphi}
&\phin\to\phi&&\mbox{weakly star in } \,W^{1,\infty}(0,T;\VD)\cap H^1(0,T;V)\cap L^\infty(0,T;W)\,,\\
\label{conmu}
&\mun\to\mu&&\mbox{weakly star in }\,L^\infty(0,T;V)\,,\\
\label{conw}
&\wn\to w&&\mbox{weakly in }\,H^1(0,T;H)\,.
\end{align}
Since, owing to the compactness of the embedding $\,W\subset C^0(\overline\Omega)$, it follows from \cite[Sect.~8,~Cor.~4]{Simon} that
$H^1(0,T;V)\cap L^\infty(0,T;W)$ is compactly embedded in $\,C^0(\overline Q)$, we may also assume that
\begin{equation}
\label{conny1}
\phin\to \phi\quad\mbox{strongly in }\,C^0(\overline Q)\,,
\end{equation}
whence, by the local Lipschitz continuity of $\,f'$,
\begin{equation}
\label{conny2}
f'(\phin)\to f'(\phi) \quad\mbox{strongly in }\,C^0(\overline Q)\,.
\end{equation}

With these strong convergence properties at hand, it follows from a standard argument (which needs no repetition here) that $(\phi,\mu,w)$
is in fact a solution to the state system \State. Moreover, we can infer from the semicontinuity properties of norms and the 
estimates shown above that there is a constant $K_0>0$, which depends only on $\norma{u}_{\L2H}$ and the data of the state system, such that 
\begin{align*}
&\|\phi\|_{ W^{1,\infty}(\VD) \cap \H1V \cap \L\infty W\cap C^0(\overline Q) } \,+\,\|\mu\|_{L^\infty(0,T;V)}  \,
+ \, \| w \|_{ \H1 H }       \,\le\,K_0\,.  
\end{align*}  
In addition, we conclude from \eqref{ss1} and elliptic regularity that $\mu\in L^2(0,T;H^3(\Omega)\cap W)$ and  
\begin{equation*}
\|\mu\|_{L^2(0,T;H^3(\Omega)\cap W)}\,\le\,C\,.
\end{equation*}
Besides, taking the time derivative in \eqref{ss2}, we can infer from comparison that also $\,\mu\in H^1(0,T;\VD)\,$
and 
\begin{equation*}
\|\mu\|_{H^1(0,T;\VD)}\,\le\,C\,.
\end{equation*}

With this, the existence of a solution $(\phi,\mu,w)$ and of a constant $K_1>0$ with the asserted properties is shown. It remains to prove 
its uniqueness. This will be done below in Theorem~{\ref{Teo1}} in the continuous dependence estimate.
\end{proof}

An immediate consequence of the uniform bound established for $\,\|\phi\|_{C^0(\overline Q)}\,$ and of the uniqueness
still to be proved below is the following.
\begin{corollary}
\label{Cor1}
Assume that {\bf (A1)}--{\bf (A3)} are fulfilled. Then there is a constant $\,K_2>0\,$ depending only on the data of the
system and $\,R\,$ such that 
\Beq
\label{ssbound2}
\max_{0\le i\le 5}\, \Bigl(\max_{j=1,2}\,\|f_j^{(i)}(\phi)\|_{C^0(\overline Q)} +
\|f^{(i)}(\phi)\|_{C^0(\overline Q)}  \Bigr)\,\le\,K_2\,,
\Eeq
whenever $(\phi,\mu,w)$ is the solution to the state system \State\ in the sense of 
Theorem~{\ref{Exist}} associated with some $\,u\in{\cal U}_R$.
\end{corollary}
\Brem
\label{Rem-J}
It is worth noting that for the proof of Theorem~\ref{Exist} it was not necessary to assume that $\,f_1(r)\,$ has a sufficiently strong
(e.g., at least quadratic) growth as $|r|\to +\infty$. Such an assumption has been made in many papers  
dealing  with regular potentials. The reason for this is that 
in our approach we avoid to test \eqref{ss2n} by $\dt\phi_n$ before sufficiently strong estimates for $\phi$ (here derived in the first estimate)
are available to handle the term \,\,$-\iint_{Q_t} f_2'(\phin)\,\dt\phin\,\,$ that arises on the \rhs. It is well possible
that a corresponding line of argumentation works also in many other cases, thus avoiding the growth assumption.
\Erem
\Brem
\label{Rem1}
Our approach in the proof of Theorem~\ref{Exist} consisted in approximating all equations including \eqref{ss3}, despite the fact 
that  we had the explicit solution~\eqref{varconst} (in terms of $u$) at our disposal. 
On the other hand, the approximation turns out to be a convenient approach 
especially if \eqref{ss3} would be replaced by a more complicated and possibly coupled PDE, still with the
 control on the \rhs. In particular, we point out  that in the paper~\cite{CoSpTr}, concerned with the viscous Cahn--Hilliard equation, 
the coefficient $\gamma$ in the analog of~\eqref{ss3} is allowed to depend on the space variable $x\in \Omega$, 
being however bounded from below by a positive constant. This setting can be considered also here without major modifications: 
of course, then the approximation~\eqref{ss3n} of~\eqref{ss3} is no longer valid and should be replaced by 
\[
\iO \gamma \, \dt w_n \, v + (\wn,v)=(u,v)\quad \mbox{ for all $\,v\in V_n\,$, a.e. in $(0,T),$}
\]
with the consequence that the resulting ODE system~\eqref{ODE3} changes into 
\[
\sum_{j=1}^n \Bigl(\iO \gamma e_j e_k \Bigr) \frac{dw_{nj}}{dt} +w_{nk}=(u,e_k) \quad  \mbox{ a.e. in }\,(0,T), \quad k=1, \ldots, n.
\]
Nonetheless, note that $\bigl(\iO \gamma e_j e_k \bigr)_{\!j,k} $, $\, j,k= 1, \ldots,n$, are the coefficients of a symmetric and positive definite 
(and thus invertible) matrix, so that the resulting modified ODE system~\eqref{ODE1}--\eqref{ODE4} is still easily solvable with a time-dependent
maximal solution.
\Erem

\subsection{An auxiliary lemma}
In this section, we show the following preparatory lemma which will prove useful in numerous estimations in the following.
\begin{lemma}
\label{Lem1}
Suppose that functions 
\begin{equation}
\label{regcoeff}
a\in 
L^2(0,T;W^{1,4}(\Omega)), \quad g\in L^2(0,T;V), \quad h\in L^2(0,T;H)\,
\end{equation} 
are given.
Then there is a unique triple \,$(\phi,\mu,w)\,$ such that
\begin{align}
\label{auxphi}
&\phi\in H^1(0,T;\VD)\cap L^\infty(0,T;V)\cap L^2(0,T;W),\\
\label{auxmu}
&\mu\in L^\infty(0,T;V),\\
\label{auxw}
&w\in H^1(0,T;H),
\end{align}
as well as
\begin{align} 
\label{aux1}
&\langle \dt\phi,v\rangle\,+\iO \nabla\mu\cdot\nabla v=0 &&\mbox{for all \,$v\in V$ and a.e. in }\,(0,T)\,,\\
\label{aux2}
&-\Delta\phi-\mu-w=a\,\phi+g &&\mbox{a.e. in }\, Q\,,\\
\label{aux3}
&\gamma\dt w+w=h&&\mbox{a.e. in }\,Q\,,\\
\label{aux4}
&\dn\phi=0&&\mbox{a.e. on }\,\Sigma\,,\\
\label{aux5}
&\phi(0)=0,\quad w(0)=0,&&\mbox{a.e. in }\,\Omega\,.
\end{align}
\Accorpa\Auxiliary aux1 aux5 
Moreover, there is some constant $K_3>0$, which increases monotonically with respect to the value of the norm
$\,\|a\|_{
L^2(0,T;W^{1,4}(\Omega))}$,   
such that
\begin{align}
\label{auxbound}
&\|\phi\|_{H^1(0,T;\VD)\cap L^\infty(0,T;V)\cap L^2(0,T;W)}\,+\,\|\mu\|_{L^2(0,T;V)}\,+\,\|w\|_{H^1(0,T;H)}\non\\
&\le\,K_3 \Big(\|g\|_{L^2(0,T;V)}\,+\,\|h\|_{L^2(0,T;H)}\Big)\,.
\end{align}
\end{lemma}
\Bdim
The existence proof is again performed via a Faedo--Galerkin approximation using the same finite-dimensional spaces as in the
proof of Theorem~{\ref{Exist}}. For the sake of brevity, we avoid writing the system explicitly here and 
just provide the relevant formal a priori estimate \eqref{auxbound} for the 
continuous system, which in the rigorous argument has to be performed for the finite-dimensional approximations. 
With this a priori estimate \eqref{auxbound} at hand, the standard limit process using weak and weak-star compactness
can be carried out to prove the existence of a solution $(\phi,\mu,w)$ having the regularity \eqref{auxphi}--\eqref{auxw}.
Notice that also the uniqueness of the solution immediately follows: indeed, if $(\phi_i,\mu_i,w_i)$,
$i=1,2$, are two solutions, then $(\phi,\mu,w):=(\phi_1-\phi_2,\mu_1-\mu_2,w_1-w_2)$ satisfies the system \eqref{aux1}--\eqref{aux5}
with $g=h=0$, and \eqref{auxbound} yields that $\phi=\mu=w=0$.

To begin with, we
first note that \eqref{aux1} implies that $\overline{\dt\phi}=0$ a.e. in $(0,T)$ which, thanks to the initial condition
$\,\phi(0)=0$, yields that $\,\overline{\phi(t)}=0\,$ for all $t\in [0,T]$. We thus may test \eqref{aux1} by 
$\,\frac 12 \CN(\dt\phi)+\frac 12 \mu+\phi$, \eqref{aux2} by $\,\dt\phi-\Delta\phi-\mu$, and
\eqref{aux3} by $\,K\,\dt w$, where the constant $K>0$ is yet to be determined. Addition of the resulting equations and 
integration over $(0,t)$, where $t\in (0,T]$, then leads to the cancellation of some terms, and upon rearraging the terms, we arrive at the identity
\begin{align}
\label{eaux1}
&\frac 12 \iot\|\dt\phi(s)\|_*^2\,ds \,+\,\frac 12\,\|\phi(t)\|_V ^2\,+\,\frac K2\,\|w(t)\|_H^2\,+\,\frac 12\iint_{Q_t}|\nabla\mu|^2\non\\
&\quad+ \,\iint_{Q_t} \bigl(|\Delta\phi|^2+|\mu|^2+K\gamma|\dt w|^2\bigr)\non\\
&=\,\iint_{Q_t}\nabla\phi\cdot\nabla\mu\,+\iint_{Q_t} w\,\dt\phi\,-\iint_{Q_t} w\Delta\phi\,-\iint_{Q_t} \mu w\,
+\iint_{Q_t} a\phi(-\Delta\phi-\mu)\non\\
&\qquad + \iint_{Q_t} a\phi\dt\phi \,+\iint_{Q_t} g(-\Delta\phi-\mu)\,+\iint_{Q_t} g\dt\phi\,+\,K\iint_{Q_t} h\dt w 
:=\,\sum_{j=1}^9 I_j\,,
\end{align}
with obvious meaning. Five of the terms on the \rhs\ can be easily estimated using Young's inequality \eqref{Young}. Namely, we have
\begin{align}
\separa
\label{eaux2}
|I_1|\,&\le\,\frac 14\iint_{Q_t}|\nabla\mu|^2\,+\iint_{Q_t}|\nabla\phi|^2\,,
\\[1mm]
\separa
\label{eaux3}
|I_3|\,&\le\,\frac 14\iint_{Q_t}|\Delta\phi|^2\,+\iint_{Q_t}|w|^2\,,
\\[1mm]
\separa
\label{eaux4}
|I_4|\,&\le\,\frac 14\iint_{Q_t}|\mu|^2\,+\,\iint_{Q_t}|w|^2\,,
\\[1mm]
\separa
\label{eaux5}
|I_7|\,&\le\,\frac 14\iint_{Q_t}(|\Delta\phi|^2+|\mu|^2)\,+\,2\iint_{Q_t}|g|^2\,,
\\[1mm]
\separa
\label{eaux6}
|I_9|\,&\le\,\frac{K\gamma}4\iint_{Q_t}|\dt w|^2\,+\,\frac K{\gamma}\iint_{Q_t}|h|^2\,.
\end{align} 
The remaining four terms, which involve $\,\dt\phi\,$ and/or $\,a$, require more attention. At first, we use \eqref{umba1} and 
Young's inequality to see that
\begin{equation}
\label{eaux7}
|I_8|\,\le\,C\iot \|\dt\phi(s)\|_*\,\|g(s)\|_V\,ds\,\le\,\frac 18\iot \|\dt\phi(s)\|_*^2\,ds\,+\,C\,\|g\|^2_{L^2(0,t;V)}\,.
\end{equation} 

Next, we have
$\,\,I_2=\iO w(t)\phi(t)-\iint_{Q_t} \phi\,\dt w$, so that, by Young's inequality \eqref{Young},
\begin{equation}
\label{eaux8}
|I_2|\,\le\,\frac 14 \|\phi(t)\|_H^2\,+\,\|w(t)\|_H^2\,+\,\frac {K\gamma}4\iint_{Q_t}|\dt w|^2\,+\,\frac K\gamma\iint_{Q_t}|\phi|^2\,. 
\end{equation} 

Now observe that we have the continuous embedding $W^{1,4}(\Omega)\subset L^\infty(\Omega)$. Hence, 
$$
\|a\|_{L^2(0,T;\Linfty)}\,\le\,C\,\|a\|_{L^2(0,T;W^{1,4}(\Omega))}\,.
$$
Therefore, by Young's inequality,
\begin{align}
\label{eaux9}
|I_5|\,&\le\,\iot\|a(s)\|_\infty\,\|\phi(s)\|_H\,\bigl(\|\Delta\phi(s)\|_H\,+\,\|\mu(s)\|_H\bigr)\,ds \non\\
&\le\,\frac 14 \iint_{Q_t}(|\Delta\phi|^2+|\mu|^2)\,+\,C\iot\|a(s)\|^2_{W^{1,4}(\Omega)}\,\|\phi(s)\|_H^2\,ds\,.
\end{align}

Finally, we estimate $\,I_6$. Using \eqref{umba1}, as well as the Young and H\"older inequalities, we infer that
\begin{align}
\label{eaux10}
|I_6|\,&\le\,C\iot\|\dt\phi(s)\|_* \,\|a(s)\phi(s)\|_V\,ds\non\\
&\le\,\frac 18\iot\|\dt\phi(s)\|_*^2\,ds\,+\,C\iint_{Q_t}|a|^2\,(|\phi|^2+|\nabla\phi|^2)\,+\,C\iint_{Q_t}|\nabla a|^2\,|\phi|^2\non\\
&\le\,\frac 18\iot\|\dt\phi(s)\|_*^2\,ds\,+\,C\iot\|a(s)\|_\infty^2\,\|\phi(s)\|_V^2\,ds
\,+\,C\iot\|\nabla a(s)\|_4^2\,\|\phi(s)\|_4^2\,ds\non\\
&\le\,\frac 18\iot\|\dt\phi(s)\|_*^2\,ds\,+\,C\iot\|a(s)\|_{W^{1,4}(\Omega)}^2\,\|\phi(s)\|_V^2\,ds\,.
\end{align}

At this point, we choose $\,K=4$. It then follows from \eqref{eaux1}--\eqref{eaux10}, that there are constants
$\,C_1>0$, $C_2>0$, which do not depend on $\,a,$ $\,h\,$ and $\,g\,$, such that
\begin{align*}
&\|\phi(t)\|_V^2\,+\,\|w(t)\|_H^2\,+\iot\|\dt\phi(s)\|_*^2\,ds\,+\iint_{Q_t}\bigl(|\Delta\phi|^2+|\mu|^2+|\nabla\mu|^2+|\dt w|^2\bigr)\\
&\le\,C_1\iot \bigl(\|g(s)\|_V^2+\|h(s)\|_H^2\bigr)\,ds \,+\,C_2\iot\bigl(1\,+\,\|a(s)\|^2_{W^{1,4}(\Omega)}\bigr)
\bigl(\|\phi(s)\|_V^2+\|w(s)\|_H^2\bigr)\,ds\,.
\end{align*}
Since the real-valued function $\,\,s\mapsto C_2(1+\|a(s)\|^2_{W^{1,4}(\Omega)})\,\,$ belongs to $L^1(0,T)$, we can apply 
Gronwall's lemma, whence the inequality \eqref{auxbound} follows. In addition, the standard form of the Gronwall inequality ensures that
 the constant $\,K_3 \,$ can be chosen to be monotonically increasing with respect to $\|a\|_{L^2(0,T;W^{1,4}(\Omega))}$.
\Edim
\Brem
\label{Rem-Pier}
We point out that the assumption $a \in L^2(0,T;W^{1,4}(\Omega))$ in \eqref{regcoeff} is set for convenience, to be used in the following, 
but it can be replaced by the more general assumption $a \in L^2(0,T;W^{1,p}(\Omega))$ with $p>3$. Indeed, the above estimates in the proof 
can be repeated without major changes. In particular, the estimate \eqref{eaux10} of $I_6$ can be arranged as follows:
\begin{align*}
|I_6|\,&\le \,\frac 18\iot\|\dt\phi(s)\|_*^2\,ds\,+\,C\iot\|a(s)\|_\infty^2\,\|\phi(s)\|_V^2\,ds
\,+\,C\iot\|\nabla a(s)\|_p^2\,\|\phi(s)\|_{q}^2\,ds\\
&\le\,\frac 18\iot\|\dt\phi(s)\|_*^2\,ds\,+\,C\iot\|a(s)\|_{W^{1,p}(\Omega)}^2\,\|\phi(s)\|_V^2\,ds\,,
\end{align*}
as $ p$ is greater than the space dimension $3$ and $q : = 2p/(p-2) < 6 $, so that $V\subset L^q(\Omega)$ with continuous embedding.
\Erem

\subsection{Continuous dependence and uniqueness }
Next, we state a continuous dependence result that, in particular, guarantees the uniqueness of the solution provided by Theorem~{\ref{Exist}}. 

\Bthm
\label{Teo1}
Suppose that the conditions {\bf (A1)}--{\bf (A3)} are fulfilled. Then there exists a constant $K_4>0$ such that the 
following holds true: whenever
 $u_i\in {\cal U}_R $, $i=1,2$, are given  
and $(\phi_i,\mu_i,w_i)$, $i=1,2$, are corresponding solutions to \State\ in the sense of 
Theorem~{\ref{Exist}}, then
\begin{align}
  &\|\phi_1-\phi_2\|_{H^1(0,T;\VD)\cap \C0V \cap\L2W}
   \,+\,\|\mu_1-\mu_2\|_{L^2(0,T;V)}\,+\, \norma{w_1-w_2}_{\H1H}\non\\
	&
  \leq K_4\, \norma{u_1-u_2}_{\L2H}.
  \label{contdep1}
	\end{align}
\Ethm

\Bdim
Let us set, for convenience,
\Bsist
  && u := u_1 - u_2 \,, \quad
 \phi := \phi_1 - \phi_2 \,, \quad
  \mu := \mu_1 - \mu_2 \,, \quad
  w := w_1 - w_2\,.
  \non
\Esist
Then \,$\phi (0) = 0$ \,and\, $w(0) = 0$ \,a.e. in $\,\Omega$, as well as $\,\dn\phi=0\,$ a.e. on $\Sigma$. In addition,  
writing \eqref{ss1}--\eqref{ss3} for $(\phi_i,\mu_i,w_i)$, $i=1,2$, and taking the differences, we obtain that
\begin{align}
\label{diff1}
& \langle \partial_t \phi\, v \rangle +\int_{\Omega}\nabla \mu \cdot\nabla v =0 &&\hbox{for every $v\in V$\, and a.e. in $\,(0,T)$\,,}
\\
\label{diff2} 
& -\Delta\phi-\mu-w=  -(f'(\phi_1) -  f'( \phi_2) ) &&\mbox{a.e. in }Q\,, 
\\[1mm]
\label{diff3} 
&\gamma \, \dt w \,+\,w = u&&\mbox{a.e. in }Q\,.
\end{align} 
\quad Now observe that 
$$-(f'(\phi_1)-f'(\phi_2))\,=-\int_0^1\frac d{ds}f'(\phi_2+s(\phi_1-\phi_2))\,ds \,=\,a\,\phi\,,$$
where 
\begin{equation}
\label{diff4}
a:=-\int_0^1 f''(\phi_2+s(\phi_1-\phi_2))\,ds\,.
\end{equation}
With this choice of $\,a$, we see that the triple $(\phi,\mu,w)$ satisfies a system of the form \eqref{aux1}--\eqref{aux5} with
$\,g=0\,$ and $\,h=u$. By virtue of Lemma~{\ref{Lem1}}, the assertion will thus be proved if we can show that there exists some constant 
$\,C>0$, which depends only on the data of the system and $R$, such that
\begin{align}
\label{paolouno}
\|a\|_{L^2(0,T;W^{1,4}(\Omega))}\,\le\,C.
\end{align}
\quad Now recall  that $\,u_1,u_2\in\CU_R$. Since the constant $K_1$ from \eqref{ssbound1} depends for controls belonging to 
$\CU_R$ only on the data and $R$, it follows that  
$$\|\phi_2+s(\phi_1-\phi_2)\|_{C^0(\overline Q)}\,\le K_1\,\quad\mbox{for all }\,s\in[0,1].$$
By the continuity of $f''$, it then follows that \,$\|a\|_{C^0(\overline Q)}\,$ is bounded by a constant that only 
depends on the data and $R$. The same then holds for $\,\|a\|_{L^2(0,T;L^4(\Omega))}$. 
Finally, we obviously have that
$\,\,|\nabla a|\,\le\,C\,(|\nabla \phi_1|+|\nabla\phi_2|)\,\,$ a.e. in $Q$. But this implies that 
$\,\|\nabla a\|_{L^\infty(0,T;L^6(\Omega)^3)}\,$ is bounded, which then also holds for $\|a\|_{L^2(0,T;W^{1,4}(\Omega))}$. 

With this, the assertion is proved: note that the space $\L\infty V$ in \eqref{auxbound} is replaced by $\C0V$ in \eqref{contdep1} since 
$\phi_1, \phi_2$ are known to be continuous from $[0,T]$ to $V$ (cf.~\eqref{ssbound1}). In particular, from \eqref{contdep1}, in the case when $u_1=u_2$, it follows that $\phi=\mu=w=0$, which proves the uniqueness of the solution.
\Edim


\section{Differentiability of the control-to-state operator}
\setcounter{equation}{0}

Let us introduce the Banach spaces
\begin{align}
\label{defX}
\CX\,&:=\left(H^1(0,T;\VD)\cap L^\infty(0,T;V)\cap L^2(0,T;W)\right)\times L^2(0,T;V)\times H^1(0,T;H)\,,\\[1mm]
\label{defY}
\CY\,&:=\left(W^{1,\infty}(0,T;\VD) \cap H^1(0,T;V) \cap L^\infty(0,T; W) \cap C^0(\overline Q)\right) \nonumber\\
&\quad\,\,\, \times \left(H^1(0,T;\VD)\cap L^\infty(0,T;V)\cap L^2(0,T;W\cap H^3(\Omega))\right) \times \H1 H\,.    
\end{align}
From Theorem~{\ref{Exist}} and Theorem~{\ref{Teo1}} we know that the control-to-state operator
$$\CS:u\mapsto \CS(u)=(\CS_1(u),\CS_2(u),\CS_3(u)):=(\phi,\mu,w)\,$$ is well defined as a mapping from 
$\,\CU=L^2(0,T;H)$ into $\,\CY\,$ and Lipschitz continuous
as a  mapping from  $\,\CU_R\,$ into \,$\CX$\, for every $R>0$. 
In this section, we study the differentiability properties of this operator.
More precisely, we want to show that under the assumptions {\bf (A1)}--{\bf (A3)}  
the operator $\CS$ is  twice continuously Fr\'echet differentiable on $\CU$ 
as a mapping from $\CU$ into $\CX$. We first show the following result.
\Bthm
\label{TeoFre1}
Suppose that the conditions {\bf (A1)}--{\bf (A3)} are fulfilled. 
Then the  control-to-state operator
$\,\CS\,$ is for any $R>0$ Fr\'echet differentiable in $\CU_R$ as a mapping from $\CU$ into $\CX$.  
Moreover, for every $\us\in\CU_R$ and every increment $\,h\in\L2H$, the 
triple $(\xi,\eta,v)=\CS'(\us)[h]\in\CX$ is the unique solution to the linearized system 
\begin{align}
\label{ls1}
&\langle \dt \xi,v\rangle \,+\iO \nabla\eta\cdot\nabla v = 0 &&\mbox{for all $\,v\in V\,$ and a.e. in }\,(0,T),\\
\label{ls2}
& -\Delta \xi - \eta - v = - f''(\phis)\xi &&\mbox{a.e. in }\,Q,\\
\label{ls3}
&\gamma \dt v + v = h  &&\mbox{a.e. in }\,Q,\\
\label{ls4}
&\dn\xi = 0  &&\mbox{a.e. on }\,\Sigma,\\
\label{ls5}
&\xi(0)=0,\,\quad v(0)=0 &&\mbox{a.e. in }\,\Omega.
\end{align}
\Accorpa\Linear ls1 ls5
\Ethm
\Bdim
The existence of a unique solution $(\xi,\eta,v)\in\CX$ to the system \Linear\ follows directly from 
Lemma~{\ref{Lem1}}: indeed, the system \Linear\
is of the form \eqref{aux1}--\eqref{aux5} with $\,g=0$\, and $\,a=-f''(\phis)$, and, in view of \eqref{ssbound1}, it is easily verified that
$$
\|-f''(\phis)\|_{L^2(0,T;W^{1,4}(\Omega))}\, \le\,C,
$$
with a constant $C>0$ which depends only on the data of the state system and $\,R$. Moreover, it follows from \eqref{auxbound}
that the linear mapping $\, h\mapsto (\xi,\eta,v)\,$ is continuous from $\CU$ into $\CX$. 	

To show the asserted Fr\'echet differentiability, we consider increments $\,h\in\CU\,$ with $\,\us+h\in \CU_R\,$ and denote by $C>0$ constants
that may depend on the data and $R$, but not on the special choice of such increments. We also set
$\,(\phi^h,\mu^h,w^h):=\CS(\us+h)$, \,\,$(\phis,\mus,\ws):=\CS(\us)$, and
$$\yh:=\phi^h-\phis-\xi,\quad \zh:=\mu^h-\mus-\eta, \quad \zeh:=w^h-\ws-v.$$
We then have to show that
\begin{align}
\label{oh-oh}
&\|\CS(\us+h)-\CS(\us)-\CS'(\us)[h]\|_\CX\,=\,\|(\yh,\zh,\zeh)\|_\CX\,=\,o\bigl(\|h\|_{L^2(0,T;H)}\bigr)\non \\
&\quad\mbox{as }\,\|h\|_{L^2(0,T;H)}\,\to\,0\,.
\end{align}
Observe that \eqref{contdep1} implies that
\begin{align}
\label{F11}
&\|\phi^h-\phis\|_{H^1(0,T;\VD)\cap L^\infty(0,T;V)\cap L^2(0,T;W)}\,+\,\|\mu^h-\mus\|_{L^2(0,T;V)}\,+\,\|w^h-\ws\|_{H^1(0,T;H)}\non\\
&\le\,C\,\|h\|_{\L2H}\,.
\end{align} 
Moreover, the triple $(\yh,\zh,\zeh)\in\CX$ is obviously a solution to the system
\begin{align}
\label{fre1}
&\langle \dt\yh,v\rangle\,+\iO\nabla\zh\cdot\nabla v=0&&\mbox{for all $v\in V$ and a.e. in}\,(0,T),\\[1mm]
\label{fre2}
&-\Delta\yh-\zh-\zeh=-\bigl(f'(\phi^h)-f'(\phis)-f''(\phis)\xi\bigr)&&\mbox{a.e. in }\,Q,\\[1mm]
\label{fre3}
&\gamma\dt \zeh+\zeh=0&&\mbox{a.e. in }\,Q,\\[1mm]
\label{fre4}
&\dn\yh=0&&\mbox{a.e. on }\,\Sigma,\\[1mm]
\label{fre5}
&\yh(0)=0,\quad \zeh(0)=0,&&\mbox{a.e. in }\,\Omega,
\end{align}
whence it immediately follows that $\zeh =0$ a.e. in $Q$. Moreover, we infer from Taylor's theorem with integral remainder that
\begin{equation}
\label{Taylor1}
f'(\phi^h)-f'(\phis)-f''(\phis)\xi=f''(\phis)\yh + A^h (\phi^h-\phis)^2\quad\mbox{a.e. in}\, Q,
\end{equation}
with the remainder
\begin{equation}
\label{Ah1}
A^h:=\int_0^1 (1-s)\,f'''(\phis+s(\phi^h-\phis))\,ds\,.
\end{equation}
From this we conclude that the system \eqref{fre1}--\eqref{fre5} is of the form \eqref{aux1}--\eqref{aux5} with $\,\,a:=-f''(\phis)\,\,$
and \,\,$g:=-A^h(\phi^h-\phis)^2$.
 
In view of \eqref{ssbound1}, we have \,\,$\|a\|_{L^2(0,T;W^{1,4}(\Omega))}\,=\,\|-f''(\phis)\|_{L^2(0,T;W^{1,4}(\Omega))}
\,\le\,C$. It thus follows from \eqref{auxbound} in Lemma~{\ref{Lem1}} that \eqref{oh-oh}, and thus the assertion the theorem, is valid provided
we can show that 
\begin{equation}
\label{F12}
\|g\|_{L^2(0,T;V)}^2\,=\,\|-\,A^h(\phi^h-\phis)^2\|_{L^2(0,T;V)}^2\,\le\,C\,\|h\|^4_{L^2(0,T;H)}\,.
\end{equation} 
Now observe that, a.e. in $Q$,
\begin{align}
\label{Ah2}
&|A^h|\,\le\,C\,, \quad\,|\nabla A^h|\,\le\,C\left(|\nabla \phis|+|\nabla \phi^h|\right)\,,\\[1mm]
\label{Ah3}
&\nabla \left(A^h(\phi^h-\phis)^2\right)\,=\,\nabla A^h(\phi^h-\phis)^2\,+\,2A^h(\phi^h-\phis)\nabla(\phi^h-\phis)\,.
\end{align}
Therefore, we have that
\begin{align*}
&\|-\,A^h(\phi^h-\phis)^2\|_{L^2(0,T;V)}^2\\
&\le\,C\iint_Q |\phi^h-\phis|^4\,+\,
C\iint_Q\bigl(|\nabla \phis|^2+|\nabla\phi^h|^2\bigr)
\,|\phi^h-\phis|^4\\
&\quad\ {}+C\iint_Q |\phi^h-\phis|^2\,|\nabla(\phi^h-\phis)|^2\,=: J_1+J_2+J_3,
\end{align*}
with obvious meaning. Now, owing to \eqref{ssbound1}, \eqref{F11}, H\"older's inequality, and the continuity of the embeddings
$\,V\subset L^6(\Omega)\subset L^4(\Omega)$,
$$
J_1\,\le\,C\ioT\|(\phi^h-\phis)(t)\|_4^4\,dt\,\le\,C\,\|\phi^h-\phis\|_{L^\infty(0,T;V)}^4\,\le\,C\,\|h\|^4_{L^2(0,T;H)}\,,
$$
as well as 
\begin{align*}
J_3\,&\le\,C\ioT \|(\phi^h-\phis)(t)\|_4^2\,\|\nabla(\phi^h-\phis)(t)\|^2_{L^4(\Omega)^3}\,dt \\
&\le\,C\,\|\phi^h-\phis\|^2_{L^\infty(0,T;V)}\,\|\phi^h-\phis\|^2_{L^2(0,T;W)}\,\le\,C\,\|h\|^4_{L^2(0,T;H)}\,.
\end{align*}
Finally, we infer that
\begin{align*}
J_2\,&\le\,C\ioT \bigl(\|\nabla\phis(t)\|^2_{L^6(\Omega)^3}\,+\,\|\nabla\phi^h(t)\|^2_{L^6(\Omega)^3}\bigr)
\,\|(\phi^h-\phis)(t)\|_6^4\,dt\\
&\le\,C\,\|\phi^h-\phis\|_{L^\infty(0,T;V)}^4\,\left(\|\phis\|^2_{L^2(0,T;W)}\,+\,\|\phi^h\|^2_{L^2(0,T;W)}\right)\\
&\le\,C\,\|h\|^4_{L^2(0,T;H)}\,,
\end{align*}
which concludes the proof of the assertion.
\Edim

As the next step, we show that the mapping $\,\CS':L^2(0,T;H)\to\CL(L^2(0,T;H);\CX)$, $u\mapsto\CS'(u)$, is locally 
Lipschitz continuous. We have the following result.
\Bthm
\label{Teo2}
Suppose  that {\bf (A1)}--{\bf (A3)} are fulfilled. Then there is a constant $K_5>0$, which depends only on the data of the 
state system and on $R$, such that the following holds: whenever $\,u_i\in\CU_R$, $i=1,2$, are given, then it holds for every 
$\,h\in\L2H\,$ that
\begin{equation}
\label{contdep3}
\|(\CS'(u_1)-\CS'(u_2))[h]\|_\CX\,\le\,K_5\,\|u_1-u_2\|_{\L2H}\,\|h\|_{\L2H}\,.
\end{equation} 
\Ethm
\Bdim
Let  $h\in \L2H$ be fixed. We set $\,(\phi_i,\mu_i,w_i):=\CS(u_i)\,$ and $\,(\xi_i,\eta_i,v_i):=\CS'(u_i)[h]$, for $i=1,2$, and put
\begin{align*}
&u:=u_1-u_2,\quad \phi:=\phi_1-\phi_2, \quad \mu:=\mu_1-\mu_2,\quad w:=w_1-w_2,\quad\\
&\xi:=\xi_1-\xi_2,\quad \eta:=\eta_1-\eta_2,\quad v:=v_1-v_2\,.
\end{align*}
It then easily follows that $(\xi,\eta,v)\in\CX$ is a solution to the system
\begin{align}
\label{kurt1}
&\langle \dt\xi,\rho \rangle\,+\iO\nabla\xi\cdot\nabla \rho=0&&\mbox{for all $\rho\in V$ and a.e. }\,t\in(0,T)\,,\\[1mm]
\label{kurt2}
&-\Delta\xi-\eta-v=a\xi+g&&\mbox{a.e. in }\,Q\,,\\
\label{kurt3}
&\gamma\dt v+v=0&&\mbox{a.e. in }\,Q\,,\\
\label{kurt4}
&\dn\xi=0&&\mbox{a.e. on }\,\Sigma\,,\\
\label{kurt5}
&\xi(0)=0,\quad v(0)=0,&&\mbox{a.e in }\,\Omega\,,
\end{align} 
where we have put
\begin{equation}
a:=-f''(\phi_1), \quad g:=-(f''(\phi_1)-f''(\phi_2))\xi_2\,.
\end{equation}
Again, this is a system of the form \eqref{aux1}--\eqref{aux5}, and once more it is easily shown that $\,\,\|a\|_{L^2(0,T;W^{1,4}(\Omega))}\,\,$
is bounded. Hence, by Lemma~{\ref{Lem1}}, the result will be proved once we can show that 
\begin{equation}
\label{kurt6}
\|g\|_{\L2V}\,\le\,C\,\|u\|_{\L2H}\,\|h\|_{\L2H}\,.
\end{equation}
Now observe that, by Taylor's formula,
$$
f''(\phi_1)-f''(\phi_2)\,=\,\int_0^1 f'''(\phi_2+s(\phi_1-\phi_2))\,ds\,\phi \,=:\,B^h\,\phi\,,
$$
which in view of \eqref{ssbound2} implies that, a.e. in $Q$,
\begin{align*}
|g|\,&\le\,|B^h|\,|\xi_2|\,|\phi|\,\le\,C\,|\xi_2|\,|\phi|\,,\\
|\nabla g|\,&\le\,C\left(|\xi_2|\,|\phi|\,\bigl(|\nabla\phi_1|+|\nabla\phi_2|\bigr)\,+\,|\xi_2|\,|\nabla\phi|\,+\,|\nabla\xi_2|\,|\phi|\right)\,.
\end{align*}
Next, we recall that 
$$\|(\xi_i,\eta_i,v_i)\|_{\CX}\,=\,\|\CS'(u_i)[h]\|_{\CX}\,\le\,C\,\|h\|_{\L2H}, \quad\mbox{for \,$i=1,2$}.
$$
We therefore can conclude as follows:
\begin{align*}
\|g\|^2_{\L2H}\,&\le\,C\iQ|\xi_2|^2\,|\phi|^2\,\le\,C\ioT\|\xi_2(t)\|_4^2\,\|\phi(t)\|_4^2\,dt\non\\
&\le \,C\,\|\xi_2\|^2_{L^\infty(0,T;V)}\,\|\phi\|^2_{L^\infty(0,T;V)}\,\le\,C\,\|u\|^2_{\L2H}\,\|h\|^2_{\L2H}\,,
\end{align*}
where we also have used \eqref{contdep1}. Moreover, by similar reasoning, and using the embedding $V\subset L^4(\Omega)$ once more,
\begin{align*}
\separa
\iQ|\nabla g|^2\,&\le\,C\iQ|\xi_2|^2\,|\phi|^2\bigl(|\nabla\phi_1|^2+|\nabla\phi_2|^2\bigr)\,+\,C\iQ |\xi_2|^2\,|\nabla\phi|^2\\
&\quad \,+\,C\iQ|\nabla\xi_2|^2\,|\phi|^2\\
\separa
&\le\,C\,\|\xi_2\|_{L^\infty(0,T;V)}^2\,\|\phi\|^2_{L^\infty(0,T;V)}\,\bigl(\|\phi_1\|^2_{L^2(0,T;W)}\,+\,\|\phi_2\|^2_{\L2W}\bigr)\\
\separa
&\quad\,+\,C\,\|\xi_2\|_{L^\infty(0,T;V)}^2\,\|\phi\|^2_{L^2(0,T;W)}\,+\,C\,\|\phi\|_{L^\infty(0,T;V)}^2\,\|\xi_2\|^2_{L^2(0,T;W)}\\[1mm]
\separa
&\le\,C\,\|\xi_2\|_{L^\infty(0,T;V)\cap{\L2W}}^2\,\|\phi\|^2_{L^\infty(0,T;V)\cap \L2W}\\[1mm]
&\le\,C\,\|u\|_{\L2H}^2\,\|h\|_{\L2H}^2\,.
\end{align*}
The assertion is thus proved.
\Edim
 
Having shown this continuous dependence estimate, we can now proceed to prove that the control-to-state operator has a second
Fr\'echet derivative. We have the following result:
\Bthm
\label{TeoFre2}
Suppose that the conditions {\bf (A1)}--{\bf (A3)} are fulfilled. 
Then the  control-to-state operator
$\,\CS\,$ is for any $R>0$ twice Fr\'echet differentiable in $\CU_R$ as a mapping from $\CU$ into $\CX$.  
Moreover, for every $\us\in\CU_R$ and all increments $\,h,k\in\L2H$, the 
triple $(\psi,\nu,z)=\CS''(\us)[h,k]\in\CX$ is the unique solution to the bilinearized system 
\begin{align}
\label{bilin1}
&\langle \dt \psi,v\rangle +\iO \nabla\nu\cdot\nabla v = 0 &&\mbox{for all \,$v\in V\,$ and a.e. in }\,(0,T),\\
\label{bilin2}
&-\Delta \psi  - \nu - z = - f''(\phis)\psi - f'''(\phis)\xi^h\xi^k  &&\mbox{a.e. in }\,Q,\\
\label{bilin3}
&\gamma \dt z + z = 0  &&\mbox{a.e. in }\,Q,\\
\label{bilin4}
&\dn\psi = 0 &&\mbox{a.e. on }\,\Sigma,\\
\label{bilin5}
&\psi(0)=0,\,\quad z(0)=0 && \mbox{a.e. in }\,\Omega,
\end{align}
\Accorpa\Bilinear bilin1 bilin5
where $\,(\xi^h,\eta^h, v^h ):=\CS'(\us)[h]\,$ and $\,(\xi^k,\eta^k, v^k ):=\CS'(\us)[k]$. 
\Ethm
\Bdim
By virtue of Lemma~\ref{Lem1}, we first establish the existence of a unique solution $(\psi,\nu,z)\in\CX$ to the system \Bilinear, where we 
immediately note that $z= 0$ a.e.\ in $Q$, due to \eqref{bilin3} and \eqref{bilin5}. Indeed, the system
\Bilinear\ is of the form \Auxiliary, where in this case we have that $\,\,a:=-f''(\phis)$\, and $\,g:=-f'''(\phis)\xih\xik$. Since,
again, \,$\|a\|_{L^2(0,T;W^{1,4}(\Omega))}\,\le\,C$, it suffices to show that \,\,$\|g\|_{\L2V}\,\le\,C$. We achieve this  
by proving an estimate of the form
\begin{equation}
\label{bilin6}
\|g\|_{\L2V}\,\le\,C\,\|h\|_{\L2H}\,\|k\|_{\L2H},
\end{equation}
which in view of \eqref{auxbound} then also implies that the mapping $\,(h,k)\mapsto (\psi,\nu, z)\,$ is continuous from 
$\L2H\times\L2H$ into $\CX$. 
To this end, recall that
$$\|\xih\|_{H^1(0,T;\VD)\cap L^\infty(0,T;V)\cap \L2W}\,\le\,\|\CS'(u)[h]\|_{\CX}\,\le\,C\,\|h\|_{\L2H}\,,$$
and a corresponding estimate holds true for $\,\xik$. Therefore, using \eqref{ssbound2}, we have that
\begin{align*}
\|g\|^2_{\L2H}\,&\le\,C\iQ |\xih|^2\,|\xik|^2\,\le\,C\ioT\|\xih(t)\|_4^2\,\|\xik(t)\|^2_4\,dt\\
&\le\,C\,\|\xih\|^2_{L^\infty(0,T;V)}\,\|\xik\|^2_{L^\infty(0,T;V)}\,\le\,C\,\|h\|^2_{\L2H}\,\|k\|^2_{\L2H}\,.
\end{align*}
Moreover, in view of \eqref{ssbound2}, we have a.e. in $\,Q\,$ that
$$|\nabla g|\,\le\,C\,\left(|\nabla\phis|\,|\xih|\,|\xik|\,+\,|\nabla\xih|\,|\xik|\,+\,|\xih|\,|\nabla\xik|\right),$$
so that 
\begin{align*}
\iQ|\nabla g|^2\,&\le\,C\iQ\left(|\nabla\phis|^2\,|\xih|^2\,|\xik|^2\,+\,|\nabla\xih|^2\,|\xik|^2\,+\,|\xih|^2\,|\nabla\xik|^2\right)\\
&\le\,C\ioT\|\nabla\phis(t)\|^2_{L^6(\Omega)^3}\,\|\xih(t)\|_6^2\,\|\xik(t)\|_6^2\,dt\\
&\quad\,+\,C\ioT\left(\|\nabla\xih(t)\|_{L^4(\Omega)^3}^2\,\|\xik(t)\|_4^2\,+\,\|\xih(t)\|_4^2\,\|\nabla\xik(t)\|_{L^4(\Omega)^3}^2\right)\,dt.
\end{align*}
Hence, it results that 
\begin{align*}
\iQ|\nabla g|^2&\le\,C\,\|\phis\|^2_{L^\infty(0,T;W)}\,\|\xih\|_{L^\infty(0,T;V)}^2\,\|\xik\|^2_{L^\infty(0,T;V)}\\[1mm]
&\quad\,+\,C\left(\|\xih\|^2_{\L2W}\,\|\xik\|^2_{L^\infty(0,T;V)}\,+\,\|\xih\|_{L^\infty(0,T;V)}^2\,\|\xik\|^2_{\L2W}\right)\\[1mm]
&\le\,C\,\|h\|_{\L2H}^2\,\|k\|^2_{\L2H}\,,
\end{align*}
which concludes the existence and uniqueness proof.

We now show the differentiability result. For this, we have to show that
\begin{align}
\label{ohoh2}
&\sup_{\|h\|_{\L2H}=1}\,\|\CS'(\us+k)[h]-\CS'(\us)[h]-(\psi,\nu,z)\|_{\CX}\,=\,o\left(\|k\|_{\L2H}\right)\non\\
&\quad\mbox{as }\,\|k\|_{\L2H}\to 0\,.
\end{align}

To this end, let $\,h,k\in\L2H$\, be given with $\|h\|_{\L2H}=1$ and $\,\us+k\in\CU_R$. Next, we put 
$\,(\xih,\eta^h,\vh):=\CS'(\us)[h]\,$ and \,$(\oxih,\oetah,\ovh):=\CS'(\us+k)[h]\,.$
We have, since $\|h\|_{\L2H}=1$,
\begin{equation}
\label{F21}
\|(\xih,\eta^h,\vh)\|_{\CX}\,+\,\|(\oxih,\oetah,\ovh)\|_{\CX}\,\le\,C\,.
\end{equation}
Moreover, it follows from \eqref{contdep3} that
\begin{equation}
\label{F22}
\|(\oxih-\xih,\oetah-\eta^h,\ovh-\vh)\|_{\CX}\,\le\,C\,\|k\|_{L^2(0,T;H)}\,.
\end{equation}

Next, we consider the functions
$$\Phi:=\oxih-\xih-\psi,\quad \rho:=\oetah-\eta^h-\nu, \quad \omega:=\ovh-\vh-z\,.$$
A little calculation then shows that the triple $\,(\Phi,\rho,\omega)\in\CX\,$ solves the system
\begin{align}
\label{Fre21}
&\langle\dt\Phi,v\rangle\,+\iO\nabla\Phi\cdot\nabla v=0&&\mbox{for all $\,v\in V$\, and a.e. }\,t\in(0,T)\,,\\[1mm]
\label{Fre22}
&-\Delta\Phi-\rho-\omega\,=\,-f''(\phis)\Phi+g&&\mbox{a.e. in }\,Q\,,\\[1mm]
\label{Fre23}
&\gamma\dt \omega+\omega=0&&\mbox{a.e. in }\,Q\,,\\[1mm]
\label{Fre24}
&\dn\Phi=0&&\mbox{a.e. on }\,\Sigma\,,\\[1mm]
\label{Fre25}
&\Phi(0)=0,\quad \omega(0)=0,&&\mbox{a.e. in }\,\Omega\,,
\end{align} 
where
\begin{align}
\label{Fre26}
g\,&=\,-\bigl(f''(\phi^k)-f''(\phis)\bigr)\bigl(\oxih-\xih\bigr)\,-\,\bigl(f''(\phi^k)-f''(\phis)-f'''(\phis)\xik\bigr)\xih\non\\
&=:\,g_1+g_2\,,
\end{align}
with obvious notation. Clearly, \eqref{Fre21}--\eqref{Fre25} is again of the form \Auxiliary, and since $\,\|-f''(\phis)\|_{L^2(0,T;W^{1,4}
(\Omega))}\,\le\,C$, the assertion will be proved once we  can show that
\begin{equation}
\label{Fre27}
\|g\|^2_{\L2V}\,\le\,C\,\|k\|^4_{\L2H}\,.
\end{equation}
At first, similar estimates as above, using \eqref{ssbound1}, \eqref{ssbound2}, \eqref{F21} and \eqref{F22}, yield that
\begin{align}
&\|g_1\|^2_{\L2V}\non\\
&\le\,C\iint_Q  |f''(\phi^k)-f''(\phis)|^2\,|\oxih-\xih|^2 \non\\
&\quad\,+\,C\iQ\left(|f'''(\phi^k)|^2|\nabla(\phi^k-\phis)|^2\,+\,|\nabla\phis|^2\,|f'''(\phi^k)-f'''(\phis)|^2 \right)\,|\oxih-\xih|^2\non\\
&\quad\,+\,C\iQ|f''(\phi^k)-f''(\phis)|^2\,|\nabla(\oxih-\xih)|^2,\non
\end{align}
which leads to
\begin{align}
\label{Fre28}
&\|g_1\|^2_{\L2V}\non\\
&\le\,C\ioT \|(\phi^k-\phis)(t)\|_4^2\,\|(\oxih-\xih)(t)\|^2_4\,dt\non\\
&\quad\,+\,C\ioT\|\nabla(\phi^k-\phis)(t)\|^2_{L^4(\Omega)^3}\,\|(\oxih-\xih)(t)\|_4^2\,dt\non\\
&\quad\,+\,C\ioT\|\nabla\phis(t)\|_{L^6(\Omega)^3}^2\,\|(\phi^k-\phis)(t)\|_6^2\,\|(\oxih-\xih)(t)\|_6^2\,dt\non\\
&\quad\,+\,C\ioT\|(\phi^k-\phis)(t)\|_4^2\,\|\nabla(\oxih-\xih)(t)\|^2_{L^4(\Omega)^3}\,dt\non\\
&\le\,C\,\|k\|^4_{\L2H}\,.
\end{align}

Next, observe that
$$g_2\,=\,-\xih\left(f'''(\phis)\,(\phi^k-\phis-\xik) + Q^k\,(\phi^k-\phis)^2\right)$$
where
$$Q^k:=\int_0^1(1-s)f^{(4)}(\phis+s(\phi^k-\phis))\,ds$$
satisfies, by virtue of \eqref{ssbound2},
\begin{align}
\label{Fre29}
|Q^k|\,\le\,C, \quad |\nabla Q^k|\,\le \,C\,(|\nabla\phis|+|\nabla\phi^k|), \quad\mbox{a.e. in }\,Q.
\end{align}
We therefore have
\begin{align}
\|g_2\|^2_{\L2V}\,&\le\,C\iQ|\xih|^2\,\left(|\phi^k-\phis-\xi^k|^2\,+\,|\phi^k-\phis|^4\right)\non\\
&\quad\,+\,C\iQ|\nabla\xih|^2\,\left(|\phi^k-\phis-\xi^k|^2\,+\,|\phi^k-\phis|^4\right)\non\\
&\quad\,+\,C\iQ|\xih|^2\,\left(|\nabla\phis|^2\,|\phi^k-\phis-\xi^k|^2\,+\,|\nabla(\phi^k-\phis-\xik)|^2\right)\non\\
&\quad\,+\,C\iQ  |\xih|^2\,\left(|\nabla Q^k|^2\,|\phi^k-\phis|^4\,+\,|Q^k|^2\,|\phi^k-\phis|^2\,|\nabla(\phi^k-\phis)|^2\right).\non
\end{align}
Based on this, we can infer that
\begin{align}
\label{Fre210}
&\|g_2\|^2_{\L2V}\non \\
&\le\,C\ioT\|\xih(t)\|_4^2\,\|(\phi^k-\phis-\xi^k)(t)\|_4^2\,dt
\,+\,C\ioT \|\xih(t)\|_6^2\,\|(\phi^k-\phis)(t)\|_6^4\,dt\non\\
&\quad\,+\,C\ioT\|\nabla\xih(t)\|_{L^4(\Omega)^3}^2\,\|(\phi^k-\phis-\xik)(t)\|_4^2\,dt\non\\
&\quad\,+\,C\ioT\|\nabla\xih(t)\|_{L^6(\Omega)^3}^2\,\|(\phi^k-\phis)(t)\|_6^4\,dt\non\\
&\quad\,+\,C\ioT \|\xih(t)\|_6^2\,\|\nabla\phis(t)\|^2_{L^6(\Omega)^3}\,\|(\phi^k-\phis-\xik)(t)\|^2_6\,dt\non\\
&\quad\,+\,C\ioT\|\xih(t)\|_4^2\,\|\nabla(\phi^k-\phis-\xik)(t)\|^2_{L^4(\Omega)^3}\,dt\non\\
&\quad\,+\,C\ioT \|\xih(t)\|_6^2\,\|\nabla Q^k(t)\|^2_{L^6(\Omega)^3}\,\|(\phi^k-\phis)(t)\|_6^2\,\|(\phi^k-\phis)(t)\|_\infty^2\,dt\non\\
&\quad\,+\,C\ioT\|\xih(t)\|_6^2\,\|(\phi^k-\phis)(t)\|_6^2\,\|\nabla(\phi^k-\phis)(t)\|_6^2\,dt\non\\
&=:\,\sum_{j=1}^8 \,M_j\,,
\end{align}
with obvious notation. It remains to show that \,$\,M_j\,\le\, C\,\|k\|^4_{\L2H}$, for $1\le j\le 8$. In order not to overload the exposition,
we restrict ourselves to show this for
only two of the terms, leaving the check of the others to the interested reader. To this end, recall that in the proof of 
Theorem~{\ref{TeoFre1}}
we have shown (with $\,k\,$ replaced by $\,h\,$) that
\begin{align}
\label{Fre211}
&\|\phi^k-\phis-\xik\|_{H^1(0,T;\VD)\cap L^\infty(0,T;V)\cap \L2W}\non\\[1mm]
&\le\,\|\CS(\us+k)-\CS(\us)-\CS'(\us)[k]\|_{\CX}\,\le\,C\,\|k\|^2_{\L2H}\,.
\end{align}
By virtue of \eqref{F21} and the continuity of the embedding $V\subset L^4(\Omega)$, we therefore conclude that
\begin{align*}
|M_6|\,&\le\,C\,\|\xih\|_{L^\infty(0,T;V)}^2\,\|\phi^k-\phis-\xik\|^2_{\L2W}\,\le\,C\,\|k\|^4_{\L2H}\,.
\end{align*} 
Moreover, invoking \eqref{F21}, \eqref{Fre29}, \eqref{ssbound1}, \eqref{contdep1}, and the continutity of the embeddings 
$V\subset L^6(\Omega)$ and $W\subset L^\infty(\Omega)$, we also have that
\begin{align*}
|M_7|\,&\le\,C\,\left(\|\phis\|^2_{L^\infty(0,T;W)}\,+\|\phi^k\|^2_{L^\infty(0,T;W)}\right)
\|\phi^k-\phis\|^2_{L^\infty(0,T;V)}\,\|\phi^k-\phis\|_{\L2W}^2\\
&\le\,C\,\|k\|^4_{\L2H}\,.
\end{align*}
With this, the assertion is proved.
\Edim

Finally, we show that the mapping $\,u\mapsto \CS''(u)$\, is locally Lipschitz continuous. We have the following result. 
\Bthm 
\label{Teo3}
The mapping $\,\CS'':\L2H \to  {\cal L}(L^2(0,T;H),{\cal L}(\L2 H,\CX))$, $u \mapsto \CS''(u)$, is 
Lipschitz continuous in the following sense: there exists a constant
$K_6>0$, which depends only on $R$ and the data, such that, for all controls $u_1,u_2 \in\CU_R$ and all 
increments $h,k \in\L2 H$, it holds that 
\begin{align}
\label{contdep4}
&\|\left( \CS''  (u_1)- \CS'' (u_2)\right)[h,k]\|_{\CX}  
\le\,K_6\,\|u_1-u_2\|_{\L2H}\,\|h\|_{\L2H}\,\|k\|_{\L2H}\,.
\end{align}
\Ethm

\Bdim
Let $u_1,u_2\in\CU_R$ and $h,k\in\L2H$ be given. We put
\begin{align*}
&(\phi_i,\mu_i,w_i):=\CS(u_i), \quad (\xi_i^h,\eta_i^h,v_i^h):=\CS'(u_i)[h], \quad (\xi_i^k,\eta_i^k,v_i^k):=\CS'(u_i)[k],\\
&(\psi_i,\nu_i,z_i):=\CS''(u_i)[h,k], \quad\mbox{for \,$i=1,2$},
\end{align*}
where we recall \eqref{contdep1}, \eqref{contdep3} and the fact that $\,\,\|\CS'(u_i)[h]\|_{\CX}\,\le\,C\,\|h\|_{\L2H}$, $i=1,2$, and that an analogous
estimate holds true for $\,\CS'(u_i)[k]$. Moreover, a little calculation shows that the triple $(\psi,\nu,z)\in\CX$ solves the system 
\begin{align}
\label{F41}
&\langle \dt\psi,y\rangle\,+\iO\nabla\psi\cdot\nabla y=0&&\mbox{for all $\,y\in V$\, and a.e. in }\,Q\,,\\[1mm]
\separa
\label{F42}
&\-\Delta\psi-\nu-z=a\psi+g&&\mbox{a.e. in }\,Q\,,\\
\separa
\label{F43}
&\gamma\dt z+z=0&&\mbox{a.e. in }\,Q\,,\\
\separa
\label{F44}
&\dn\psi=0&&\mbox{a.e. on }\,\Sigma\,,\\
\separa
\label{F45}
&\psi(0)=0,\quad z(0)=0,&&\mbox{a.e. in }\,\Omega\,,
\end{align}
where $\,a:=-f''(\phi_2)\,$ is bounded in $\,L^2(0,T;W^{1,4}(\Omega))\,$ and  
\begin{align}
\label{F46}
g &:=\,-(f''(\phi_1)-f''(\phi_2))\psi_1\,-(f'''(\phi_1)-f'''(\phi_2))\xi_1^h\,\xi_1^k\non\\
&\quad\ \ {}-f'''(\phi_2)\,(\xi_1^h-\xi_2^h)\,\xi_1^k\,-\,f'''(\phi_2)\,\xi_2^h(\xi_1^k-\xi_2^k)\non\\
&=:\,\sum_{j=1}^4 g_j\,,
\end{align}
with obvious notation. In view of \eqref{auxbound} in Lemma~{\ref{Lem1}}, it suffices to show that 
\begin{align}
\label{F47}
\|g_j\|^2_{\L2V}\,\le\,C\,\|u\|_{\L2H}\,\|h\|_{\L2H}\,\|k\|_{\L2H} \quad\mbox{for }\,j=1,2,3,4.
\end{align}  

We demonstrate this only for the second and third terms. The other two terms can be treated similarly and, in order to keep the
paper at a reasonable length, are left to the reader. We have, using \eqref{ssbound1}, \eqref{ssbound2}, \eqref{contdep1},
\eqref{contdep3}, H\"older's inequality, and the continuity of the embeddings $V\subset L^6(\Omega)\subset L^4(\Omega)$ and
$W\subset \Linfty$,
\begin{align*}
&\|g_2\|^2_{\L2V}\,\\
&\le\,C\ioT\|(\phi_1-\phi_2)(t)\|_6^2\,\|\xi_1^h(t)\|_6^2\,\|\xi_1^k(t)\|_6^2\,dt \\
&\quad\,\,+\,C\ioT \|(\phi_1-\phi_2)(t)\|^2_\infty\,\|\nabla\phi_1(t)\|_{L^6(\Omega)^3}^2 \|\xi_1^h(t)\|_6^2\,\|\xi_1^k(t)\|_6^2\,dt\\
&\quad\,\,+\,C\ioT \|\nabla(\phi_1-\phi_2)(t)\|^2_{L^6(\Omega)^3} \|\xi_1^h(t)\|_6^2\,\|\xi_1^k(t)\|_6^2\,dt\\
&\quad\,\,+\,C\ioT \|(\phi_1-\phi_2)(t)\|^2_6\,\|\nabla\xi_1^h(t)\|_{L^6(\Omega)}^2 \,\|\xi_1^k(t)\|_6^2\,dt\\
&\quad\,\,+\,C\ioT \|(\phi_1-\phi_2)(t)\|^2_6\,\|\xi_1^h(t)\|_6^2\,\|\nabla\xi_1^k(t)\|_{L^6(\Omega)^3}^2\,dt,
\end{align*}
so that
\begin{align*}
&\|g_2\|^2_{\L2V}\,\\
&\le\,C\,\|\phi_1-\phi_2\|^2_{L^\infty(0,T;V)}\,\|\xi_1^h\|^2_{L^\infty(0,T;V)}\,\|\xi_1^k\|^2_{L^\infty(0,T;V)}\\[1mm]
\separa
&\quad\,\,+\,C\,\|\phi_1\|^2_{L^\infty(0,T;W)}\,\|\xi_1^h\|^2_{\L\infty V}\,\|\xi_1^k\|^2_{\L\infty V}\,\|\phi_1-\phi_2\|^2_{\L2W}\\[1mm]
\separa
&\quad\,\,+\,C\,\|\xi_1^h\|^2_{\L\infty V}\,\|\xi_1^k\|^2_{\L\infty V}\,\|\phi_1-\phi_2\|^2_{\L2W}\\[1mm]
\separa                                                                 
&\quad\,\,+\,C\,\|\phi_1-\phi_2\|^2_{\L\infty V}\,\|\xi_1^h\|^2_{\L2W}\,\|\xi_1^k\|^2_{\L\infty V}\\[1mm]
\separa
&\quad\,\,+\,C\,\|\phi_1-\phi_2\|^2_{\L\infty V}\,\|\xi_1^h\|^2_{\L\infty V}\,\|\xi_1^k\|^2_{\L2W}
\\[2mm]
\separa
&\le\,C\,\|\CS(u_1)-\CS(u_2)\|_{\CX}^2\,\|\CS'(u_1)[h]\|_{\CX}^2\,\|\CS'(u_1)[k]\|^2_{\CX} \\[2mm]
\separa
&\le\,C\,\|u_1-u_2\|^2_{\L2H}\,\|h\|^2_{\L2H}\,\|k\|^2_{\L2H}\,.
\end{align*}
Similarly, it holds that 
\begin{align*}
\separa
\|g_3\|^2_{\L2V}\,&\le\,C\ioT\|(\xi_1^h-\xi_2^h)(t)\|_4^2\,\|\xi_1^k(t)\|_4^2\,dt\non\\
\separa
&\quad\,+\,C\ioT \|\nabla\phi_2\|^2_{L^6(\Omega)^3}\,\|(\xi_1^h-\xi_2^h)(t)\|_6^2\,\|\xi_1^k(t)\|_6^2\,dt\non\\
\separa
&\quad\,+\,C\ioT \|\nabla(\xi_1^h-\xi_2^h)(t)\|^2_{L^4(\Omega)^3}\,\|\xi_1^k(t)\|_4^2\,dt\non\\
\separa
&\quad\,+\,C\ioT \|(\xi_1^h-\xi_2^h)(t)\|_4^2\,\|\nabla\xi_1^k\|_4^2\,dt \\
&\,\le\,C\,\|u_1-u_2\|^2_{\L2H}\,\|h\|^2_{\L2H}\,\|k\|_{\L2H}^2\,.
\end{align*}
With this, the assertion is proved.
\Edim
\Brem
\label{Rem2}
With Theorem~{\ref{Teo3}}, we have shown that the control-to-state operator $\CS$ is twice continuously Fr\'echet differentiable as a 
mapping from $\CU=\L2H$ into $\CX$. This result paves the way to prove first-order necessary and second-order sufficient optimality
conditions for the optimal control problem {\bf (CP)} in the following section.  
\Erem

\section{The optimal control problem}
 
\setcounter{equation}{0}

In this section, we study the optimal control problem {\bf (CP)} with the cost functional \eqref{cost}. Besides the general
conditions {\bf (A1)}--{\bf (A3)}, we  make the following 
assumptions:
\begin{description}
\item[(A4)] \,\,It holds $b_1\ge 0$, $b_2\ge 0$, $b_3>0$, and $\,\kappa>0$. 
\item[(A5)] \,\,The thresholds $\,\underline u, \overline u\in\LiQ\,$ satisfy $\,\underline u\le\overline u\,$ almost everywhere in
$\,Q$, and the target functions satisfy $\phi_Q\in L^2(Q)$ and $\phi_\Omega \in V.$
\end{description}
We assume $\kappa > 0$ to include the effects of sparsity. By an obvious modification, the theory of second-order 
conditions remains valid also for $\kappa = 0$.
\Brem
\label{Rem3}
The assumption $\phi_\Omega \in V$ is  useful in order to have more regular solutions to the associated
adjoint system (see below). It is not overly restrictive in view of the continuous embedding \,$(H^1(0,T;H)\cap L^2(0,T;W))
\subset C^0([0,T];V)$\, which implies that $\phi(T)\in V$. 
\Erem

The following existence result can be shown with a standard argument.
\Bthm
\label{Teo4}
Let {\bf (A1)}--{\bf (A5)} hold and suppose that $G:\L2H\to\erre$ is 
nonnegative, convex and continuous. Then the optimal control problem {\bf (CP)} admits a solution 
$\us\in\Uad$.
\Ethm

\subsection{The adjoint system}

In the following, let $\us\in\CU_R$ be fixed and $(\phis,\mus,\ws)=\CS(\us)$ be the associated state.
 The corresponding adjoint state system is formally given by:
\begin{align}
\label{adj1}
&-\dt p -\Delta q + f''(\phis )q =b_1(\phis -\phi_Q) &&\mbox{in }\,Q,\\
\label{adj2}
&-\Delta p-q=0 &&\mbox{in }\,Q,\\
\label{adj3}
&-\gamma\dt r+r-q=0 &&\mbox{in }\,Q,\\
\label{adj4}
&\dn p=\dn q=0 &&\mbox{on }\,\Sigma,\\
\label{adj5}
&p(T)=b_2(\phis(T)-\phi_\Omega), \quad r(T)=0 &&\mbox{in }\,\Omega.
\end{align}
\Accorpa\Adjoint adj1  adj5
We immediately observe that the system is decoupled in the sense that $\,r\,$ can be directly recovered from \eqref{adj3} with the
terminal condition $r(T)=0$ once $q$ is determined. We point out that \eqref{adj1} has to be rewritten in a weak (variational) form. 
We now show a well-posedness result for a slightly more general system.
\Bthm
\label{Teo5}
Suppose that {\bf (A1)}--{\bf (A3)} are fulfilled, and assume that $\us\in\CU_R$ with $(\phis,\mus,\ws)=\CS(\us)$, $g_1\in\L2H$ and
$g_2\in V$ are given. Then the system 
\begin{align}
\separa
\label{wadj1} 
& \langle -\dt p,v \rangle\,+\iO\nabla q\cdot\nabla v +\iO f''(\phis)q\,v=\iO g_1\,v&&\mbox{for all \,$v\in V$, a.e. in } \,(0,T),\\
\separa
\label{wadj2}
&-\Delta p-q=0&&\mbox{a.e. in }\,Q\,,\\
\separa
\label{wadj3}
&-\gamma\dt r+r-q=0&&\mbox{a.e. in }\,Q\,,\\
\separa
\label{wadj4}
&\dn p=0&&\mbox{a.e. on }\,\Sigma\,,\\
\separa
\label{wadj5}
&p(T)=g_2,\quad r(T)=0,&&\mbox{a.e. in }\,\Omega\,,
\end{align}
has a unique solution triple $(p,q,r)$ with the regularity
\begin{align}
\separa
\label{regp}
&p\in H^1(0,T;\VD)\cap C^0([0,T];V)\cap L^2(0,T;W\cap H^3(\Omega)),\\
\separa
\label{regq}
&q\in L^2(0,T;V),\\
\separa
\label{regr}
&r\in H^1(0,T;V).
\end{align}
Moreover, there is a constant $K_7>0$, which depends only on $\,R\,$ and the data, such that the two inequalities below hold:
\begin{align} 
&\|p\|_{\C0 H \cap \L2W}\,+\,\|q\|_{L^2(0,T;H)}
\,+\,\|r\|_{H^1(0,T;H)} \non\\
&\le\, K_7\left(\|g_1\|_{L^2(0,T;H)} + \|g_2\|_H \right),\label{adj6}\\[2mm]
&\|p\|_{C^0([0,T];V)\cap L^2(0,T;W\cap H^3(\Omega))}\,+\,\|q\|_{L^2(0,T;V)}
\,+\,\|r\|_{H^1(0,T;V)} \non\\
&\le\, K_7\left(\|g_1\|_{L^2(0,T;H)} + \|g_2\|_V \right).\label{adj7}
\end{align}
\Ethm
\Bdim The linear initial-boundary value problem given by \eqref{wadj1}, \eqref{wadj2}, \eqref{wadj4}, together with the first terminal condition
in \eqref{adj5}, is again solved  via a Faedo--Galerkin approximation using the same eigenvalues, eigenfunctions and $n$-dimensional 
approximating spaces $V_n$ as in the proof of Theorem~{\ref{Exist}}. For the sake of shortness, we avoid to write the approximating $n$-dimensional
analogues of \eqref{wadj1}--\eqref{wadj2} explicitly here and just provide the relevant a priori estimates formally for the 
continuous problem. Having these estimates for the $n$-dimensional approximations, one can apply the standard weak and weak-star
compactness arguments to pass to the limit as $n\to\infty$, thereby showing the existence of the sought solution. Uniqueness then follows
immediately from the linearity and the estimate \eqref{adj6}. 

To this end, we insert $v=p$ in \eqref{wadj1} and test \eqref{wadj2} by \,\,$-q$. Then we add the resultants, noting that a 
cancellation of two terms occurs, and integrate over $(t,T)$, where $t\in[0,T)$ is arbitrary. Introducing the notation
$Q^t:=\Omega\times(t,T)$ for $t\in[0,T)$, we then obtain, after rearranging terms and invoking \eqref{ssbound2} and Young's inequality,
\begin{align*}
\separa
&\frac 12\|p(t)\|_H^2\,+\,\iint_{Q^t}|q|^2\,=\,\frac 12 \|g_2\|_H^2\,-\,\iint_{Q^t} f''(\phis)q\,p
+ \iint_{Q^t} g_1\,p \\
&\le\,\frac 12\,\left(\|g_1\|_{\L2H}^2+\|g_2\|_H^2\right) \,+\frac 12\iint_{Q^t}|q|^2\,+\,C\iint_{Q^t} |p|^2\,.
\end{align*}
Gronwall's lemma then yields that
\begin{equation}
\label{adjesti1}
\|p\|_{L^\infty(0,T;H)}^2\,+\,\|q\|_{L^2(0,T;H)}^2\,\le\,C\,\left(\|g_1\|_{\L2H}^2+\|g_2\|_H^2\right)\,.
\end{equation}
In addition, we conclude from \eqref{wadj2} and \eqref{wadj4}, invoking standard elliptic estimates, that 
\begin{equation}
\label{pier1}
\|p\|^2_{L^2(0,T;W)}\ \le\,C\,\left(\|g_1\|_{\L2H}^2+\|g_2\|_H^2\right)\,,
\end{equation}
and comparison in \eqref{wadj3} immediately shows that also
\begin{equation}
\label{pier2}
\|r\|^2_{H^1(0,T;H)}\,\le\,C\,\left(\|g_1\|_{\L2H}^2+\|g_2\|_H^2\right)\,.
\end{equation}
The validity of the inequality \eqref{adj6} is thus shown.

Next, we insert $v=q$ in \eqref{wadj1} and test \eqref{wadj2} (formally) by $\,-\dt p$. Addition and integration over $(t,T)$ then
yields
\begin{align*}
&\frac 12\iO|\nabla p (t)|^2\,+\,\iint_{Q^t}|\nabla q|^2\,=\,\frac 12\iO|\nabla g_2|^2\,+\iint_{Q^t}\left(g_1 q-f''(\phis)q^2\right)\\
&\le\, \|g_2\|_V^2 \,+\, \|g_1\|_{\L2H}^2 \, + \, C\iint_{Q^t}|q|^2\,\le\,C\left(\|g_1\|_{\L2H}^2+\|g_2\|_V^2\right)\,,
\end{align*}
by virtue of Young's inequality, and invoking the assumption $g_2\in V$ along with \eqref{ssbound2} and \eqref{adjesti1}. Then, we deduce that 
\begin{equation}
\label{adjesti2}
\|p\|_{L^\infty(0,T;V)}\,+\,\|q\|_{\L2V}\,\le\, C\,\left(\|g_1\|_{\L2H}^2+\|g_2\|_V^2\right),
\end{equation}
whence, using \eqref{wadj2} and elliptic regularity,
\begin{equation}
\label{pier3}
\|p\|_{L^2(0,T;W\cap H^3(\Omega))}\,\le\,C\left(\|g_1\|_{\L2H}^2+\|g_2\|_V^2\right).
\end{equation}
In addition,  using \eqref{adjesti2} and the endpoint condition $\,r(T)=0$, we obtain directly from \eqref{wadj3} that also
\begin{equation}
\label{adjesti4}
\|r\|^2_{H^1(0,T;V)}\,\le\,C\left(\|g_1\|_{\L2H}^2+\|g_2\|_V^2\right).
\end{equation}
Moreover, comparison in \eqref{wadj1} yields that 
\begin{equation*}
\|p\|_{H^1(0,T;\VD)}\,\le\,C\left(\|g_1\|_{\L2H}^2+\|g_2\|_V^2\right).
\end{equation*}
Finally, as the embedding $\,\bigl(H^1(0,T;\VD)\cap L^2(0,T;H^3(\Omega))\bigr)\subset C^0([0,T];V)\,$ is continuous, \eqref{regp} follows, and the inequality~\eqref{adj7} results from \eqref{adjesti2}--\eqref{adjesti4}. 
This concludes the proof of the assertion.
\Edim 
\Brem
\label{Rem4} This remark collects three different comments.\\[1mm]
1.  From the proof of Theorem~\ref{Teo5} the reader may realize that a weaker existence and uniqueness 
result holds if $g_1\in\L2H$ and $g_2$ is just in $H$, with a solution $(p,q,r)$ having
 the regularity expressed in \eqref{adj6}. In fact,
if we consider a weaker formulation of \eqref{wadj1}, namely,  
$$
\phantom{|}_{W^*}\langle -\dt p,v \rangle_{W} - \iO q\Delta v + \iO f''(\phis)q\,v=\iO g_1\,v \quad
\mbox{ for all \,$v\in W$, a.e. in } \,(0,T),
$$
then we deduce from \eqref{adjesti1}--\eqref{pier2} and a comparison of terms in the above equation that  
\begin{equation*}
\|\dt p\|_{L^2(0,T;W^*)}\,\le\,\,C\left(\|g_1\|_{\L2H}^2+\|g_2\|_H^2\right),
\end{equation*}
which, along with \eqref{pier1}, implies that $p\in C^0([0,T];H)$ and gives a meaning to the first terminal 
condition in \eqref{wadj5}. Of course, uniqueness then follows from inequality~\eqref{adj6}.\\[1mm]
2.  In the case of the adjoint system~\eqref{adj1}--\eqref{adj5}, the one of interest for our theory, we have 
$\,g_1:=b_1(\phis-\phi_Q)\in\L2H\,$ and $\,g_2:=b_2(\phis(T)-\phi_\Omega)\in V$, so that Theorem~{\ref{Teo5}} 
ensures that for every $\us\in\L2H$ there is a unique solution $\,(\ps,\qs,\rs)\,$ with the regularity
\eqref{regp}--\eqref{regr} that satisfies~\eqref{adj7}.\\[1mm]
3.  If, for $i=1,2$, $u_i\in\CU_R$ is given with the associated state
$(\phi_i,\mu_i,w_i)=\CS(u_i)$ and adjoint state $(p_i,q_i,r_i)$, then the triple $(p_1-p_2,q_1-q_2,r_1-r_2)$ solves the
system \eqref{wadj1}--\eqref{wadj5} with  $\,g_1:=b_1(\phi_1-\phi_2)\in\L2H\,$ and $\,g_2:=b_2(\phi_1(T)-\phi_2(T))\in V$. 
In view of \eqref{contdep1}, it therefore follows from \eqref{adj7} the estimate 
\begin{align}
\label{contdepadj}
&\|p_1-p_2\|_{C^0([0,T];V)\cap L^2(0,T;W\cap H^3(\Omega))} \,+\,\|q_1-q_2\|_{L^2(0,T;V)} \,+\,\|r_1-r_2\|_{H^1(0,T;V)}\nonumber\\
&\le\,C\,\|\phi_1-\phi_2\|_{C^0([0,T];V)}\,\le\,C\,\|u_1-u_2\|_{L^2(0,T;H)}\,.
\end{align}
\Erem

\subsection{First-order necessary optimality conditions}

In this section, we derive first-order necessary optimality conditions for local minima of the
optimal control problem {\bf (CP)}. We assume that {\bf (A1)}--{\bf (A5)} are fulfilled and that $G:L^2(0,T;H)\to\erre$ 
is a general nonnegative, convex and continuous functional. We define  the reduced cost functionals associated with the functionals  
$J$ and ${\cal J}$  introduced in \eqref{cost} by
\begin{equation}\label{reduced}
\widehat J(u) :=  J(\CS_1(u),u), \quad \widehat {\cal J}(u)={\cal J}(\CS_1(u),u)\,.
\end{equation}
Since $\CS=(\CS_1,\CS_2,\CS_3)$ is twice continuously Fr\'echet differentiable  from $\CU=\L2H$ into the space 
$C^0([0,T];H)\times L^2(0,T;V)\times H^1(0,T;H)$ (which contains $\CX$), 
it follows from the chain rule that the smooth part $\,\widehat J\,$ of the reduced objective functional is  a twice continuously 
 Fr\'echet differentiable mapping from $\CU$ into $\erre$, where, for every $\us\in\L2H$ and every $h\in \L2H$,  it holds 
with $(\phis,\mus,\ws)=\CS(u^*)$ that
\begin{align}
\label{DJ}
\widehat J'(\us)[h]\,&=\,b_1\iint_Q \xih(\phi^*-\phi_Q)\,+\,b_2\int_\Omega\xih(T)(\phis(T)-\phi_\Omega) 
\,+\, b_3 \iint_Q \us h \,,
\end{align}
where $\,(\xih,\eta^h,v^h)=\CS'(\us)[h]\in\CZ\,$ is the unique solution to the linearized system \Linear\
associated with $h$.

In the following, we assume that $\us\in\Uad$ is a locally optimal control for ${\bf (CP)}$ in the sense
of \,$\LiQ$. In this connection, recall that a control $\us\in\Uad$ is called \emph{locally optimal 
in the sense of $L^p(Q)$} for some $p\in[1,+\infty]$ if and only if there is some $\varepsilon>0$ such that
\begin{equation}
\label{lomin}
\widehat {\cal J}(u)\,\ge\,\widehat {\cal J}(\us)\quad\mbox{for all $u\in\Uad$ satisfying }\,\|u-\us\|_{L^p(Q)}\,\le\,\varepsilon.
\end{equation}
It is easily seen that any locally optimal control in the sense of $L^p(Q)$ with $1 \le p < \infty$
is also locally optimal in the sense of $\LiQ$. 
Therefore, a result proved for locally optimal controls in the sense of  $\LiQ$ is also valid for locally optimal controls 
in the sense of  $L^p(Q)$ for $1\le p<\infty$. In particular, it is true for (globally) optimal controls.

A standard argument (for details, see, e.g., \cite{SpTr1, SpTr2, CoSpTr}) then shows 
that there is 
some $\lambda^* \in \partial G(\us)\subset L^2(0,T;H)$ such that 
\begin{equation} \label{varineq1}
{\widehat J}'(\us)[u - \us] +  \kappa\iint_Q 
\lambda^* (u- \us)\,\ge 0 \quad \forall \,
u \in \Uad.
\end{equation}
As usual, we simplify the expression ${\widehat J}'(\us)[u-\us]$ in \eqref{varineq1}
by means of the adjoint state variables defined in \Adjoint. A standard calculation using the linearized system \Linear\
then leads to the following result.
\Bthm
\,\,{\rm (Necessary optimality condition)}  \label{Thm4.5}Suppose that {\bf (A1)}--{\bf (A5)} 
are fulfilled and that $G:L^2(0,T;H)\to\erre$ is nonnegative, convex and continuous. Moreover, 
 let $\us \in \Uad$ be a locally optimal control of  {\bf (CP)} in the sense of $\,\LiQ\,$ 
with associated state $\,(\phis,\mus,\ws)={\cal S}(\us)$
and adjoint state $\,(\ps,\qs,\rs)$.
Then there exists some $\lambda^*   \in \partial G(\us)$ such that,
for all $\, u \in \Uad$, 
\begin{align}
\label{varineq2}
&\iint_Q  \left(r^*+\kappa\lambda^*+b_3 u^*\right)\left(u-\us\right) \, \ge \,0 \,.
\end{align}
\Ethm

\subsection{Sparsity of controls}

The convex function \,$G$\, in the objective functional accounts for the sparsity of optimal controls, i.e., the possibility
that any locally optimal control may vanish in some subset of the space-time cylinder $Q$. The form of this region 
depends on the particular choice of the functional $\,G$.
The sparsity properties can be deduced from the variational inequality \eqref{varineq2} and the particular
form of the subdifferential  $\,\partial G$. Here, we restrict ourselves to the case of {\em full sparsity}
which is connected to the $L^1(Q)-$norm functional \,$G$\, introduced in \eqref{defg}. 
Its subdifferential is given by (see \cite{ioffe_tikhomirov1979})
\begin{equation}\label{dg}
\partial G(u) = \left\{\lambda \in L^2(Q):\,
\lambda(x,t) \in \left\{
\begin{array}{ll}
\{1\} & \mbox{ if } u(x,t) > 0\\
{[-1,1]}& \mbox{ if } u(x,t) = 0\\
\{-1\} & \mbox{ if } u(x,t) < 0\\
\end{array}
\right.
\mbox{ for a.e. } {(x,t) \in Q}
\right\}.
\end{equation}
With exactly the same argument as in the proof of the corresponding result \cite[Thm.~4.7]{CoSpTr}, we obtain the following result.
\Bthm 
\label{Teo6}
{\rm (Full sparsity)} \,\,\,Suppose that the assumptions {\bf (A1)}--{\bf (A5)} are fulfilled, and assume that 
$\underline u$ and $\overline u$ are constants such that $\underline u <0 <\overline u$. Let $\us\in \Uad$ be a 
locally optimal control in the sense 
of \,$\LiQ$\, for the problem {\bf (CP)}
with the functional $\,G\,$ defined in \eqref{defg}, and with associated state $(\phis,\mus,\ws)=\CS(\us)$ solving \State\ and 
adjoint state $(p^*,q^*,r^*)$ solving \Adjoint. Then there exists a function $\lambda^*\in \partial G(\us)$ 
satisfying \eqref{varineq2}, and it holds
\begin{eqnarray}
\us(x,t) = 0 \quad &\Longleftrightarrow& \quad |r^*(x,t)| \le \kappa, \quad\mbox{for a.e. }\,(x,t)\in Q. 
\label{fullsparse}  
\end{eqnarray}
Moreover, if\, $r^*$ and $\lambda^*$ are given, then
$\us$  is obtained from the projection formula
\begin{eqnarray*}
\us(x,t)& =& \max\left\{\underline u, \min\left\{\overline u, -b_3^{-1} \left(r^*+ \kappa \,\lambda^*\right)(x,t)\right\}\right\}
\,\mbox{ for a.e. $(x,t)\in Q$}.
\end{eqnarray*}
\Ethm

\subsection{Second-order sufficient optimality conditions} 
We conclude this paper with the derivation of second-order sufficient optimality conditions for functions \,$\us\,$ 
obeying the first-order necessary optimality 
conditions of Theorem~\ref{Thm4.5}. Second-order sufficient optimality conditions are based on a condition of coercivity that 
is required to hold for the smooth part  $\,J\,$ 
of $\,{\cal J}\,$ in a certain critical cone. The nonsmooth part $\,G\,$ contributes to sufficiency by its convexity. In the following,
we generally assume that the conditions {\bf (A1)}--{\bf (A5)} are fulfilled.
Our analysis will follow closely the lines of \cite{CoSpTr}, which in turn is an adaptation of the technique used in 
the proof of \cite[Thm.~3.4]{casas_ryll_troeltzsch2015} for the sparse control of the FitzHugh--Nagumo system. 

To this end, we fix  a control $\,\us \,$ satisfying the first-order necessary optimality conditions,
and we set $\,(\phis,\mus,\ws)=\CS(\us)$. Then the cone
\[
C(\us) = \{ v \in L^2(0,T;H) \,\text{ satisfying the sign conditions \eqref{sign} a.e. in $Q$}\},
\]
where
\begin{equation} \label{sign}
v(x,t) \left\{
\begin{array}{l}
\ge 0 \,\, \text{ if }\,\, u^*(x,t) = \underline u\\ 
\le 0 \,\, \text{ if }\,\, u^*(x,t) = \overline u 
\end{array}
\right. \,,
\end{equation}
is called the {\em cone of feasible directions}, which is a convex and closed subset of $L^2(0,T;H)$.
 We also need the directional derivative of $\,G\,$ at $u\in L^2(0,T;H)$ in the direction 
$v\in L^2(0,T;H)$, which is given by
\begin{equation}
\label{g'}
G'(u,v) = \lim_{t \searrow 0} \frac{1}{t}(G(u+t v)-G(u))\,.
\end{equation}
Following the definition of the critical cone in \cite[Sect.~3.1]{casas_ryll_troeltzsch2015}, we define
\begin{equation}
\label{critcone}
C_{\us} = \{v \in C(\us): \widehat{J}'(\us)[v] + \kappa G'(\us,v) = 0\}\,,
\end{equation}
which is also a closed and convex subset of $L^2(0,T;H)$. According to \cite[Sect.~3.1]{casas_ryll_troeltzsch2015}, 
it consists of all $v\in C(\us)$ satisfying  
\begin{equation} \label{pointwise}
v(x,t) \left\{
\begin{array}{l}
= 0 \,\,\mbox{ if } \,\,|r^*(x,t) + b_3 u^*(x,t)| \not= \kappa\\
\ge 0 \,\,\mbox{ if }\,\, u^*(x,t) = \underline u\,\, \mbox{ or } \,\,(r^*(x,t) = -\kappa \,\,\mbox{ and }\,\, u^*(x,t) = 0)\\
\le 0 \,\,\mbox{ if }\,\, u^*(x,t) = \overline u\,\, \mbox{ or }\,\, (r^*(x,t) = \kappa \,\,\mbox{ and } \,\,u^*(x,t) = 0)
\end{array}
\right.\,.
\end{equation}

At this point, we provide an explicit expression for $\,\widehat J''(u)[h,k]\,$ 
for arbitrary $\,u,h,k\in\L2H$. Arguing exactly as in the derivation of the corresponding formula \cite[Eq.~(4.52)]{CoSpTr}, 
we obtain that
\begin{align}	
& \widehat J''(\us)[h,k] \,=\,\iint_Q \bigl(b_1-f^{(3)}(\phis)q^*\bigr)\xih\,\xik \,+\,b_2\iO \xih(T)  \xik(T)
\,+\,b_3\iint_Q h \,k\,,
\label{walter1}
\end{align}
where $\,(\xih,\eta^h,\vh)=\CS'(\us)[h]\,$ and $\,(\xik,\eta^k,v^k)=\CS'(\us)[k]$. 

For the proof of the second-order sufficient 
optimality condition, we will need the following preparatory result.
\begin{lemma} 
\label{Lem4.7}
Assume that {\bf (A1)}--{\bf (A5)} are satisfied and that $\us\in\Uad$ is given with associated state $(\phis,\mus,\ws)$ and
adjoint state $(\ps,\qs,\rs)$. 
Suppose that $\{\widetilde u_j\}\subset\Uad$ converges strongly in $L^2(0,T;H)$ to $\,\us$, and that $\{h_j\}\subset \L2H$
converges weakly in $\L2H$ to $h$. In addition, let $(\widetilde\phi_j,\widetilde\mu_j,\widetilde w_j)=
\CS(\widetilde u_j)$, and let $\,(\widetilde p_j,\widetilde q_j,\widetilde r_j)$ be the associated adjoint state, for $j\in\enne$. Moreover, let, 
for arbitrary $h\in \L2 H$, $(\xih,\eta^h,v^h)=\CS'(\us)[h]$, as well as $(\widetilde \xi^{h_j},\widetilde\eta^{h_j},
\widetilde v^{h_j})=\CS'(\widetilde u_j)[h_j]$, for $j\in\enne$.
 Then
\begin{align} 
\label{rudi1}
&\lim_{j\to\infty}\,\widehat J'(\widetilde u_j)[h_j]\,=\,\widehat J'(\us)[h]\,,\\[2mm]
\label{rudi2}
&\lim_{j \to \infty} \Big(\iint_Q\bigl(b_1-f^{(3)}(\widetilde\phi_j)\widetilde q_j\bigr)\bigl|\widetilde\xi^{h_j}\bigr|^2
\,+\,b_2\iO \bigl|\widetilde\xi^{h_j}(T)\bigr|^2 \Big) \nonumber\\
&\quad =\iint_Q\bigl(b_1-f^{(3)}(\phis)q^*\bigr)\bigl|\xi^h\bigr|^2\,+\,b_2\iO\bigl|\xi^h(T)\bigr|^2\,.
\end{align} 
\end{lemma}
\begin{proof}
At first, notice that \eqref{contdepadj} yields that $\,\|\widetilde r_j-r^*\|_{H^1(0,T;H)}\to 0\,$ as $\,j\to\infty$, which implies that
$$
\lim_{j\to\infty}\,\widehat J'(\widetilde u_j)[h_j]\,=\,\lim_{j\to\infty}\iint_Q (\widetilde r_j+b_3 \widetilde u_j)h_j
\,=\,\iint_Q(\rs+b_3\us)h\,=\,\widehat J'(\us)[h]\,,
$$         
i.e., \eqref{rudi1} is valid. Next, in order to prove \eqref{rudi2}, we observe that
$$
(\widetilde\xi^{h_j}, \widetilde\eta^{h_j}, \widetilde v^{h_j})-(\xi^h,\eta^h,v^h) = \left(\CS'(\widetilde u_j)
-\CS'(\us)\right)[ h_j]\, +\,\CS'(\us)[h_j-h]\,.
$$
By virtue of \eqref{contdep3} and the boundedness of $\{h_j\}$ in $L^2(0,T;H)$, the 
first summand on the right converges to zero strongly in $\,\CX$.
The second converges to zero weakly star in (cf.~\eqref{defX})
$$ \bigl(H^1(0,T;\VD)\cap \L\infty V\cap L^2(0,T;W)\bigr) \times L^2(0,T;V) \times H^1(0,T;H). $$
Thanks to the compact embedding $V \subset L^p(\Omega)$ for $1\le p<6$, the compactness result stated in 
\cite[Sect.~8, Cor.~4]{Simon}) then ensures that
\begin{equation}
\label{rudi3}
\widetilde\xi^{h_j}\to \xi^h \quad\mbox{strongly in }\,C^0([0,T];L^5(\Omega))\,.
\end{equation}
In particular, we have that
\begin{align} 
\label{rudi4}
&\lim_{j\to \infty} \Big(b_1\iint_Q \bigl|\widetilde\xi^{h_j}\bigr|^2\,+\,b_2\iO\bigl|\widetilde\xi^{h_j}(T)\bigr|^2 \Big)
= b_1 \iint_Q \bigl|\xi^h\bigr|^2 \,+\,b_2\iO\bigl|\xi^h(T)\bigr|^2\,.
\end{align}
Moreover, we obtain from \eqref{contdep1} that $\,\|\widetilde \phi_j-\phis\|_{H^1(0,T;\VD)\cap L^\infty(0,T;V)}\to 0$ as $\,j\to\infty$, so 
that we can conclude from \cite[Sect.~8, Cor.~4]{Simon}, the
global estimate \eqref{ssbound2}, and \eqref{contdepadj}, that, as $j\to \infty$,
\begin{align}
\label{rudi5}
&\|f^{(3)}(\widetilde\phi_j)-f^{(3)}(\phis)\|_{C^0([0,T];L^5(\Omega))} \,\to\,0,\\
\label{rudi6}
&\|\widetilde q_j-q^*\|_{L^2(0,T;H)}\,\to\,0.
\end{align}
Combining this with \eqref{rudi3}, we readily verify that
\begin{equation}
\lim_{j\to\infty}\iint_Q f^{(3)}(\widetilde\phi_j)\widetilde q_j\bigl|\xi^{h_j}\bigr|^2\,=\,\iint_Q
f^{(3)}(\phis)q^*\bigl|\xi^h\bigr|^2,
\end{equation} 
which concludes the proof.
\end{proof}

With Lemma~{\ref{Lem4.7}} shown, the road is paved for the proof of second-order sufficient optimality conditions. To this 
end, we will employ the following coercivity condition:

\begin{equation} \label{coerc}
\widehat J''(\us)[v,v] >0 \quad \forall \, v \in C_{\us} \setminus \{0\}\,.
\end{equation}

\noindent Condition \eqref{coerc} is a direct extension of associated conditions that are standard in finite-dimensional 
nonlinear optimization. In the optimal control of partial differential equation, it was first used in \cite{casas_troeltzsch2012}.
We have the following result.
\begin{theorem}
\label{Teo7} \,\,{\rm (Second-order sufficient condition)} \,\,Suppose that {\bf (A1)}--{\bf (A5)} are fulfilled. 
Moreover, let $\us \in \Uad$, together with the associated state 
$(\phi^*,\mu^*,w^*)=\CS(\us)$ and the adjoint state $(p^*,q^*,r^*)$, fulfill the first-order necessary optimality conditions 
of Theorem~\ref{Thm4.5}. If, in addition,  $\us$ satisfies the coercivity condition \eqref{coerc}, 
then there exist constants $\,\varepsilon > 0\,$ and $\,\zeta > 0\,$ such that the quadratic growth condition
\begin{equation} \label{growth}
\widehat {\cal J}(u) \ge \widehat {\cal J}(\us) + \zeta \, \|u-\us\|^2_{\L2 H} 
\end{equation}
holds for all $\,u \in \Uad\,$ with $ \|u-\us\|_{\L2 H}  < \varepsilon$. Consequently,
$\,\us\,$ is a locally optimal control in the sense of $\,\L2 H$.
\end{theorem}
\begin{proof}
The proof is exactly the same as that of the corresponding \cite[Thm.~4.8]{CoSpTr}. In order not to seem repetitive, 
we therefore only sketch the argument, pointing out the places in the proof where the results of 
Lemma~{\ref{Lem4.7}} are used.
We argue by contradiction, assuming that the claim of the theorem is not true. Then there exists a sequence of controls 
$\{u_j\}\subset \Uad$ such that, for all $j\in\enne$, 
\begin{equation}
\label{contrary}
\|u_j-\us\|_{\L2 H} < \frac{1}{j} \quad \mbox{ while } \quad \widehat {\cal J}(u_j)<  \widehat {\cal J}(\us) 
+ \frac{1}{2j} \|u_j-\us\|_{\L2 H}^2\,.
\end{equation}
Noting that $u_j \not=\us$ for all $j\in\enne$, we define 
$$
\tau_j := \|u_j-\us\|_{\L2 H}
  \quad \mbox{ and } \quad h_j := \frac{1}{\tau_j}(u_j-\us)\,.
$$                                               
Then $\|h_j\|_{\L2 H}=1$ and, possibly after selecting a subsequence, we can assume that 
\[
h_j \to h \, \mbox{ weakly in }\,  \L2 H
\]
for some $h\in \L2 H$. The proof is now split into three parts. 

(i) \,$h \in C_{\us}$: \,\,Here, one has to show that (cf. \eqref{critcone}) \,$\widehat J'(\us)[h]+\kappa G'(\us,h)=0$,
in particular. This follows exactly as in the corresponding step (i) in \cite{CoSpTr}, where in the proof of the
inequality \,$\widehat J'(\us)[h]+\kappa G'(\us,h)\le 0$ the identity \eqref{rudi1} is used.

(ii) \,$h = 0$: \,\,The proof of this claim is again exactly the same as in the corresponding step (ii) in the proof
of \cite[Thm.~4.8]{CoSpTr}. It uses the identity \eqref{rudi2} in order to show that \eqref{contrary}
implies that $\,\widehat J''(\us)[h,h]\le 0$, whence $h=0$ follows using \eqref{coerc}.

(iii) {\em Contradiction:} 
Again, the argumentation is exactly the same as the corresponding step (iii) in the proof of
\cite[Thm.~4.8]{CoSpTr}:
we know from the previous step that $\,h_j\to 0\,$ weakly in $\,\L2H$. Now, \eqref{walter1} yields that
\begin{align}
\label{umba-umba}
&\widehat J_1''(\us)[h_j,h_j] = \iint_Q \bigl(b_1-f^{(3)}(\us)q^*\bigr)|\xi^{h_j}|^2\,+\,b_2\iO|\xi^{h_j}(T)|^2
\,+\,b_3\iint_Q |h_j|^2\,,
\end{align}
where we have set $\,(\xi^{h_j},\eta^{h_j},v^{h_j})=\CS'(\us)[h_j]$, for $\,j\in\enne$. Since $h_j\to 0$ weakly in
$\L2 H$, we find from \eqref{rudi2} in Lemma~\ref{Lem4.7} that the sum of the first two integrals on the right-hand side of 
\eqref{umba-umba} converges to zero. 
On the other hand, $\,\|h_j\|_{\L2 H}=1$ for all 
$\,j\in\enne$, by construction. Therefore, 
\begin{align}
\label{uff-uff}
&\liminf_{j\to\infty} \,\widehat J''(\us)[h_j,h_j] \,\ge \,\liminf_{j\to\infty} \,b_3\iint_Q|h_j|^2\,=\,
b_3\,>0\,.
\end{align}
On the other hand, the condition \eqref{contrary} leads, using \eqref{rudi2} again, to the conclusion that
$$\liminf_{j\to\infty} \,\widehat J''(\us)[h_j,h_j] \,\le \,0\,,$$
which contradicts \eqref{uff-uff}. The assertion of the theorem is thus proved. 
\end{proof}

\smallskip

\section*{Acknowledgments}
This research has been supported by the MIUR-PRIN Grant
2020F3NCPX ``Mathematics for industry 4.0 (Math4I4)''. 
In addition, {PC gratefully acknowledges his affiliation 
to the GNAMPA (Gruppo Nazionale per l'Analisi Matematica, 
la Probabilit\`a e le loro Applicazioni) of INdAM (Isti\-tuto 
Nazionale di Alta Matematica). 

\bigskip


\End{document}